\documentclass[]{article}
\usepackage{graphicx}
\usepackage{amsfonts}
\usepackage{amsmath}
\usepackage{amssymb}
\usepackage{fancyhdr}
\usepackage{titlesec}
\usepackage{indentfirst}
\usepackage{booktabs}
\usepackage{verbatim}
\usepackage{color}
\usepackage{amsthm}
\usepackage{caption}
\usepackage{subcaption}
\usepackage{physics}
\usepackage{mathtools}
\usepackage{tikz}
\usetikzlibrary{shapes.geometric}
\usetikzlibrary{fit}
\usetikzlibrary{calc}
\usepackage{enumerate}
\usepackage{stmaryrd}

\usepackage{hyperref}
\usepackage{cleveref}
\usepackage[page,header]{appendix}
\usepackage{titletoc}

\numberwithin{equation}{section}
\newtheorem{theorem}{Theorem}[section]

\newtheorem{remark}[theorem]{Remark}

\newcommand{\be}{\boldsymbol{\eta}}

\newcommand{\bp}{\mathbf{p}}

\newcommand{\e}{\mathrm{e}}

\newcommand{\bx}{\mathbf{x}}

\newcommand{\dt}{\Delta t}

\renewcommand{\dd}{\mathrm{d}}
\newcommand{\ii}{\mathrm{i}}
\newcommand{\fourier}{\mathcal{F}_\bp}
\newcommand{\bz}{\mathbf{z}}

\renewcommand{\vec}[1]{\operatorname{vec}\left(#1\right)}
\newcommand{\vecinv}[1]{\operatorname{vec}^{-1}\left(#1\right)}

\tikzset{
  core/.style  ={circle, draw=black, thick, fill=white, inner sep=0pt, minimum size=0.42cm},
  ltri/.style  ={isosceles triangle, isosceles triangle apex angle=75, draw=black, thick,
                 fill=white, inner sep=0pt, minimum height=0.4cm},            
  rtri/.style  ={ltri, shape border rotate=180},                              
  diam/.style  ={diamond, draw=black, fill=black, inner sep=0pt, minimum size=0.28cm},
  fbox/.style  ={rectangle, draw=black, thick, fill=white, inner sep=0pt, minimum size=0.34cm},
  bondlab/.style={font=\scriptsize, below=1pt},
  leglab/.style ={below},
}

\begin{document}

\title{Dynamical Tensor Train Approximation for Wigner--Fokker--Planck Equation with Applications to \\Open Quantum Systems}

\author{
Geshuo Wang\footnote{Department of Applied Mathematics, University of Washington, Seattle, WA 98195 (geshuo@uw.edu).} 
\ and 
Jingwei Hu\footnote{Department of Applied Mathematics, University of Washington, Seattle, WA 98195 (hujw@uw.edu).}
}

\maketitle

\begin{abstract}
We develop a dynamical tensor-train method for the Wigner--Fokker--Planck equation arising in continuous-variable open quantum systems.
The Wigner function is evolved directly in tensor-train form using a projector-splitting integrator, while the nonlocal Wigner pseudo-differential operator is treated through low-rank representations of the potential difference and matrix-free exponential actions.
The resulting method avoids full tensor-product computations in high-dimensional phase space.
Numerical experiments in one, two, and three spatial dimensions demonstrate the accuracy and scalability of the method, including simulations in a six-dimensional phase space.
\end{abstract}

\section{Introduction}

In physics, an open quantum system refers to a quantum system coupled with external degrees of freedom \cite{breuer2007theory}.
The interaction changes the quantum dynamics significantly, leading to dissipation and diffusion.
As open quantum system models are widely applied in many fields \cite{carmichael2009open,weiss2012quantum}, the numerical simulation of such dynamical systems is important in scientific computing.
For continuous-variable systems, the Wigner phase-space formulation \cite{wigner1932quantum,hillery1984distribution} describes the reduced state by a quasi-probability distribution $W(t,\bx,\bp)$ on position-momentum space.
For a particle coupled to a thermal bath, a commonly used phase-space model for the reduced dynamics is the Wigner--Fokker--Planck (WFP) equation \cite{caldeira1983path,arnold2004analysis}, whose long-time behavior and stationary states have been studied in \cite{sparber2004long,arnold2012wigner}.
The equation combines a nonlocal Wigner pseudo-differential operator, which carries the coherent quantum effects generated by the potential, with friction and diffusion terms that model the environment.
It therefore connects open quantum dynamics with kinetic and Fokker--Planck descriptions while retaining intrinsically quantum effects, and it is a natural model for decoherence \cite{zurek2003decoherence} and quantum transport in nanoscale devices \cite{weinbub2018recent}.

The numerical simulation of the WFP equation has two main difficulties.
The first is the different structures of different parts of the equation.
The Wigner pseudo-differential operator is nonlocal and is most naturally treated in Fourier space.
The transport, friction, and diffusion terms, however, are local differential operators in phase space.
The intrinsic difference of the operators makes operator splitting a standard tool for Wigner-type problems \cite{arnold1995operator,arnold1996operator,yi2026time}.
The treatment of the local terms is relatively standard.
For the Wigner pseudo-differential operator, there are two main ideas.
The first treats it exactly, either as a convolution in momentum with the Wigner kernel by numerical quadrature \cite{shao2011adaptive}, or in Fourier space, where the operator acts as a multiplication by the potential difference \cite{ringhofer1990spectral,xiong2016advective,chen2019high}.
The second replaces the operator by its Moyal expansion \cite{moyal1949quantum}, which requires computation of high-order derivatives for both the potential and the Wigner function \cite{gardner1994quantum,li2014numerical}.

The second major challenge is the curse of dimensionality.
In $d$ physical dimensions, the Wigner function depends on $2d$ phase-space variables.
With $N_x$ and $N_p$ grid points used in each spatial and momentum direction, respectively, a direct tensor-product discretization requires $\mathcal{O}(N_x^dN_p^d)$ degrees of freedom.
For $d=3$ and $N_x=N_p=64$, this is already $64^6\approx 6.9\times 10^{10}$ unknowns, about $550$ GB in double precision.
Nevertheless, resolving the full Wigner function is desirable when one is interested in phase-space structure, reduced observables, or the interplay between system dynamics and environmental dissipation.
Although several methods have been extended beyond one dimension \cite{xiong2016advective,chen2019high}, their cost still grows rapidly with the phase-space dimension.

This challenge motivates numerical methods based on low-rank representations.
For functions of two variables, a low-rank representation is generally given by the singular value decomposition.
In higher dimensions, several tensor formats are available, including the Tucker format \cite{tucker1966some}, the canonical polyadic (CP) decomposition \cite{hitchcock1927expression}, hierarchical Tucker and tree tensor networks \cite{hackbusch2012tensor}, and the tensor-train (TT) format \cite{oseledets2011tensor}.
The TT format is also called matrix product states in the field of quantum many-body physics \cite{schollwock2011density}.
Low-rank methods have been applied successfully to high-dimensional kinetic equations \cite{einkemmer2025review} and to quantum dynamics \cite{meyer1990multi,haegeman2016unifying}.
For open quantum systems, tensor-network and related compressions have proved effective for many-body lattice systems \cite{verstraete2004matrix,sun2024simulation}, for path-integral methods \cite{strathearn2018efficient,jorgensen2019exploiting}, for hierarchical equations of motion \cite{shi2018efficient,ke2023tree}, and for bath correlation functions \cite{zhan2026reducing,zhang2026low} or bath influence functionals \cite{erpenbeck2023tensor,yu2025inchworm,wang2026accelerated}.
In addition, low-rank matrices \cite{appelo2025kraus}, tensor trains \cite{delmastro2026completely}, and multiconfiguration expansions \cite{cattarius2004multidimensional} have been applied to the Lindblad equation.
These approaches compress the density matrix, the path integral, or the bath influence functional, and not the Wigner function in phase space.

For a wider range of time-dependent differential equations, low-rank methods can be broadly classified into step-and-truncate (SAT) approaches and dynamical low-rank (DLR) approximations.
SAT methods generally advance the equation over one time step in a low-rank format and then compress the result \cite{dolgov2014low,kormann2015semi,guo2022low}.
They are conceptually simple and easily made rank-adaptive, but intermediate rank growth can be costly.
A DLR method evolves the solution directly on a prescribed low-rank manifold by projecting the equation onto its tangent space \cite{koch2007dynamical,lubich2013dynamical}.
The projector-splitting integrator splits this projection into a sequence of simpler substeps.
It was introduced for matrices \cite{lubich2014projector}, and extended to tensor trains \cite{lubich2015time} and Tucker tensors \cite{lubich2018time}, and its error is robust with respect to small singular values \cite{kieri2016discretized}.
The basis-update-and-Galerkin (BUG) integrators form a second family, which updates the low-rank bases before evolving the coefficients in the updated bases \cite{ceruti2022unconventional,ceruti2022rank,ceruti2024parallel,ceruti2024robust}.

Low-rank methods have also been applied successfully to Wigner-type and Fokker--Planck-type equations in recent years.
For the Wigner equation, a dynamical low-rank algorithm was proposed in \cite{shao2026dynamical}, in which a separable representation of the difference potential, combined with truncations of the Wigner kernel, yields a separated form of the pseudo-differential operator.
For the Wigner--Poisson system, a sampling-based adaptive-rank solver was developed in one spatial dimension, where the solution was observed to have low-rank structure for moderate and large dimensionless Planck constants \cite{christlieb2026sampling,christlieb2026conservative}, and was extended to a structure-preserving solver in hierarchical Tucker format for the full 3D3V phase space \cite{christlieb2026structure}.
Low-rank methods for classical Fokker--Planck-type equations are more abundant \cite{dolgov2012fast,coughlin2022efficient,tang2024solving,nakao2025reduced,zhang2026separable}.
Numerical examples also show that in the collisional regime the solution relaxes toward a low-rank profile \cite{wang2026dynamical,wang2026implicit}. To the best of our knowledge, these two lines of research have remained separate: existing Wigner-type approaches focus on purely coherent quantum dynamics, whereas Fokker–Planck-type approaches address classical dynamics. Neither directly addresses the open quantum setting governed by the WFP equation, where coherent quantum evolution and quantum dissipative effects interact to produce dynamics distinct from those described by either equation alone. This gap motivates the investigation of low-rank structure in the WFP equation.

In this work we develop a dynamical tensor-train method for the deterministic simulation of high-dimensional open quantum dynamics described by the WFP equation.
The Wigner function is represented in TT form, with every position and momentum variable a separate tensor mode, and is evolved by the projector-splitting integrator on the TT manifold \cite{lubich2015time}.
The nonlocal pseudo-differential operator is projected onto the current low-rank manifold and advanced in time by the action of a matrix exponential, while the local transport, friction, and diffusion terms are advanced with the forward Euler scheme in the TT representation.
The method is applicable when the potential difference is represented in CP or TT form.
The computational cost scales linearly in $N_x$, quasi-linearly in $N_p$, and cubic in the tensor rank $r$.
The cost is much lower than the one in full tensor representation, making simulations in a six-dimensional phase space feasible.
We carry out numerical examples in one-, two-, and three-dimensional settings, and observe that friction and diffusion damp fine phase-space structures and reduce the effective TT ranks.
 
The remainder of the paper is organized as follows.
\Cref{sec_2Wigner_Fokker_Plank_Equation} introduces the Wigner--Fokker--Planck equation.
\Cref{sec_Tensor_Train_representation} presents the TT representation of the Wigner function and low-rank representations of the potential difference.
\Cref{sec_projector_splitting_method} develops the projector-splitting formulation and its matrix-free implementation.
\Cref{sec_discretization} describes the position, momentum, and Fourier discretizations, and \Cref{sec_computational_cost} discusses the computational complexity.
Numerical results in one, two, and three dimensions are presented in \Cref{sec_Numerical_Examples}, and \Cref{sec_conclusion} concludes the paper.

\section{Wigner--Fokker--Planck equation}
\label{sec_2Wigner_Fokker_Plank_Equation}
In this paper, we consider the Wigner--Fokker--Planck equation for a particle of mass $m$ in a $d$-dimensional space, which describes the evolution of the Wigner function $W(t,\bx,\bp)$:
\begin{equation}
\label{eq_wfp}
\partial_t W
+
\frac{\bp}{m}\cdot \nabla_\bx W
+
\Theta[V]W
=
2\gamma \nabla_\bp \cdot (\bp W)
+
D_{pp} \Delta_\bp W
+
2 D_{xp} \nabla_\bx\cdot\nabla_\bp W
+
D_{xx}\Delta_\bx W.
\end{equation}
Here, $\bx\in\mathbb{R}^d$ and $\bp\in\mathbb{R}^d$ denote the position and momentum variables, respectively, and the Wigner function $W$ is real-valued.
The term $\frac{\bp}{m}\cdot \nabla_\bx W$ describes free transport and is analogous to the transport term in classical kinetic equations.
The term $\Theta[V]W$ denotes the Wigner pseudo-differential operator associated with the potential $V(\bx)$.
The right-hand side consists of the momentum-space friction term $2\gamma \nabla_\bp \cdot(\bp W)$, the momentum diffusion term $D_{pp} \Delta_\bp W$, the cross-diffusion term $2D_{xp} \nabla_\bx\cdot\nabla_\bp W$, and the spatial diffusion term $D_{xx} \Delta_\bx W$.

For a normalized quantum state, the Wigner function satisfies
\begin{equation}
    \int_{\mathbb{R}^d}
    \int_{\mathbb{R}^d}
    W(t,\bx,\bp)\, \dd \bx\, \dd\bp = \operatorname{tr(\rho(t))}= 1,
\end{equation}
where $\rho$ is the density operator for the quantum system.

When the right-hand side is absent, i.e.,
$\gamma=D_{pp}=D_{xp}=D_{xx}=0$,
\cref{eq_wfp} reduces to the Wigner equation, which describes the evolution of the Wigner quasiprobability distribution of a closed quantum system.
With the friction and diffusion terms included, the Wigner--Fokker--Planck equation provides a phase-space model for an open quantum system evolving under the potential $V(\bx)$.
For nonnegative diffusion coefficients $D_{pp}$ and $D_{xx}$, complete positivity of the corresponding density-operator evolution is ensured when the coefficients satisfy the Lindblad condition
\begin{equation}
\label{eq_lindblad_condition}
    D_{pp}D_{xx}-D_{xp}^2
    \geq
    \frac{\hbar^2\gamma^2}{4}.
\end{equation}

In \cref{eq_wfp}, the Wigner pseudo-differential operator is defined by
\begin{equation}
\Theta[V] W(\bx,\bp)
=
-\frac{\ii}{(2\pi)^d \hbar}
\int_{\mathbb{R}^d}
\int_{\mathbb{R}^d}
\delta V(\bx,\be)
W(\bx,\bp')
\e^{\ii (\bp-\bp') \cdot\be}
\dd \bp'\,\dd \be,
\end{equation}
where
\begin{equation}
\delta V(\bx,\be)
=
V\left(\bx+\frac{\hbar}{2}\be\right)
-
V\left(\bx-\frac{\hbar}{2}\be\right).
\end{equation}
We define the Fourier transform and inverse Fourier transform with respect to the momentum variable $\bp$ by
\begin{equation}
    \fourier[W](\bx,\be)
    =
    \int_{\mathbb{R}^d}
    W(\bx,\bp)
    \e^{-\ii\bp\cdot\be}
    \dd \bp,
    \qquad
    \fourier^{-1}[\widehat{W}](\bx,\bp)
    =
    \frac{1}{(2\pi)^d}
    \int_{\mathbb{R}^d}
    \widehat{W}(\bx,\be)
    \e^{\ii\bp\cdot\be}
    \dd \be.
\end{equation}
With these definitions, the Wigner pseudo-differential term can be written as
\begin{equation}
\label{eq_pseudo_differential_operator_Fourier}
\Theta[V]W
=
-\frac{\ii}{\hbar}
\fourier^{-1}
\left[
\delta V(\bx,\be)\fourier W
\right].
\end{equation}

For a sufficiently regular potential $V(\bx)$, the Wigner pseudo-differential operator admits the formal Moyal expansion \cite{moyal1949quantum}
\begin{equation}
\label{eq_pseudo_differential_operator_Moyal}
\Theta[V]W
=
-
\sum_{n=0}^\infty
\frac{(-1)^n\hbar^{2n}}{2^{2n}}
\sum_{|\boldsymbol{\alpha}|=2n+1}
\frac{1}{\boldsymbol{\alpha}!}
\left(
\partial_\bx^{\boldsymbol{\alpha}}V
\right)
\left(
\partial_\bp^{\boldsymbol{\alpha}}W
\right),
\end{equation}
where $\boldsymbol{\alpha}$ denotes a $d$-dimensional multi-index.
The leading term of \cref{eq_pseudo_differential_operator_Moyal} is
$-\nabla_\bx V\cdot\nabla_\bp W$.
Therefore, in the classical limit $\hbar\rightarrow0$, the Wigner pseudo-differential operator reduces to the classical force term.

In our method, we use the Fourier representation \cref{eq_pseudo_differential_operator_Fourier} to evaluate the Wigner pseudo-differential operator.
Specifically, the Fourier transform of $W$ is first multiplied pointwise by the potential difference $\delta V(\bx,\be)$ in Fourier space, followed by an inverse Fourier transform in the momentum variables.

\section{Low-rank representations}
\label{sec_Tensor_Train_representation}

\subsection{Tensor train representation of the Wigner function}
In \cref{eq_wfp}, the Wigner function $W$ is a time-dependent function of $2d$ variables.
Direct numerical simulation therefore suffers from the curse of dimensionality, particularly for $d=2,3$.
Assuming $N_x$ and $N_p$ grid points in each spatial and momentum direction, respectively, the memory cost for storing the Wigner function $W$ is $\mathcal{O}(N_x^dN_p^d)$, not to mention the associated computational cost.
To alleviate this high-dimensional computational burden, low-rank tensor techniques are employed in this paper.
Specifically, we use the tensor train decomposition to represent the solution $W$ and evolve it directly in tensor train form.

Since the Wigner function depends on $d$ spatial and $d$ momentum variables, we introduce the notation $\bz=(z_1,\cdots,z_s)$ with $s=2d$, where
\begin{equation}
z_1=x_1,\quad z_2=p_1,\quad
z_3=x_2,\quad z_4=p_2,\quad\cdots,\quad z_{s-1} = x_d,\quad z_s = p_d.
\end{equation}
A tensor train representation of $W$ is given by
\begin{equation}
\label{eq_W_TTform}
W(\bz) = \sum_{a_1=1}^{r_1}
\cdots \sum_{a_{s-1}=1}^{r_{s-1}}
G_{1a_1}^{(1)}(z_1)
G_{a_1a_2}^{(2)}(z_2) \cdots
G_{a_{s-2}a_{s-1}}^{(s-1)}(z_{s-1})
G_{a_{s-1}1}^{(s)}(z_s),
\end{equation}
where $(r_1,\cdots,r_{s-1})$ denotes the TT rank.
In the study of tensor networks, a tensor train can be represented intuitively by a tensor diagram.
For the case when $d=2$, the tensor train representation for the Wigner function can be given by 
\begin{equation}
    W(\bz) = 
    \begin{tikzpicture}[baseline=(c1.base)]
      \foreach \i/\x in {1/0,2/1.8,3/3.6,4/5.4} \coordinate (c\i) at (\x,0);
      \draw[thick] (c1) -- (c4);
      \foreach \i/\l in {1/{z_1{=}x_1},2/{z_2{=}p_1},3/{z_3{=}x_2},4/{z_4{=}p_2}}
        \draw[thick] (c\i) -- ++(0,-0.8) node[leglab]{$\l$};
      \foreach \i/\l in {1/{G^{(1)}},2/{G^{(2)}},3/{G^{(3)}},4/{G^{(4)}}}
        \node[core, label=above:{$\l$}] at (c\i) {};
      \node[bondlab] at (0.9,0) {$a_1$};
      \node[bondlab] at (2.7,0) {$a_2$};
      \node[bondlab] at (4.5,0) {$a_3$};
    \end{tikzpicture}.
\end{equation}
The time dependence is omitted here for notational simplicity.
In \cref{sec_Tensor_Train_representation,sec_projector_splitting_method}, we keep the spatial and momentum variables in continuous form and their discretizations will be introduced in \cref{sec_discretization}. 

Although the Wigner function itself is real-valued, the Wigner pseudo-differential operator involves complex arithmetic.
Therefore, in the following discussions and our implementation, we allow all entries in the tensor train to be complex-valued.

\begin{remark}
In this paper, we choose to arrange the variables in an interleaved manner.
Another natural choice is to arrange the variables in the order of $(x_1,\cdots,x_d,p_1,\cdots,p_d)$.
We would like to remark that our method also works for this order. 
In general, there is no universally superior order in this setting.
\end{remark}

For the same high-dimensional tensor, the tensor-train representation is not unique because of gauge freedom between neighboring TT cores. A convenient representation can be obtained by imposing left- and right-orthogonality conditions \cite{holtz2012manifolds}.
The standard projector splitting method of tensor train \cite{lubich2015time} also utilizes these orthogonal forms. 
We then introduce the orthogonal forms that will be used in the following discussions.

For a fixed core $k$, we use a mixed canonical representation
\begin{equation}
\label{eq_canonical_form}
W(\bz) = \sum_{a_{k-1}=1}^{r_{k-1}}
\sum_{a_k=1}^{r_k}
P_{a_{k-1}}^{(<k)}(\bz_{<k})
C_{a_{k-1}a_k}^{(k)}(z_k)
Q_{a_k}^{(>k)}(\bz_{>k}),
\end{equation}
where the left and right basis functions are orthonormal:
\begin{equation}
\label{eq_orthonormal_condition}
\left\langle P_a^{(<k)}, P_b^{(<k)} \right\rangle = \delta_{ab},
\quad
\left\langle Q_a^{(>k)}, Q_b^{(>k)} \right\rangle = \delta_{ab}.
\end{equation}
Note that we do not need to construct $P^{(<k)}$ or $Q^{(>k)}$ in full tensor form if it consists of more than one physical index.
We refer the reader to \cite{lubich2015time} for more details.

After orthogonalizing the $k$th center core, the same tensor can be written as
\begin{equation}
\label{eq_bond_form}
W(\bz) = \sum_{a_k=1}^{r_k}
\sum_{b_k=1}^{r_k}
P_{a_k}^{(\leqslant k)}(\bz_{\leqslant k})
S_{a_k b_k}^{(k)}
Q_{b_k}^{(>k)}(\bz_{>k}).
\end{equation}
In \cref{eq_canonical_form,eq_bond_form}, we use the shorthand notation
\begin{equation*}
\bz_{<k} = (z_1,\cdots,z_{k-1}),
\quad
\bz_{>k} = (z_{k+1},\cdots,z_s),
\quad
\bz_{\leqslant k} = (z_1,\cdots,z_k).
\end{equation*}
The forms in \cref{eq_canonical_form,eq_bond_form} with conditions \cref{eq_orthonormal_condition} can be obtained via a sequence of singular value decompositions (SVD) or QR decompositions from both sides of the tensor train towards the center core.
This process is standard in the applications of tensor trains.
Thus, we omit the details and refer the reader to \cite{oseledets2011tensor,holtz2012manifolds,lee2018fundamental} for more information.

If we fix $k=3$ in a four-core tensor train, \cref{eq_canonical_form} can be given by the following diagram:
\begin{equation}
    \begin{tikzpicture}[baseline=(c1.base)]
  \foreach \i/\x in {1/0,2/1.8,3/3.6,4/5.4} \coordinate (c\i) at (\x,0);
  \draw[thick] (c1) -- (c4);
  \foreach \i in {1,2,3,4} \draw[thick] (c\i) -- ++(0,-0.8);
  \node[ltri] (n1) at (c1) {};
  \node[ltri] (n2) at (c2) {};
  \node[core] (n3) at (c3) {};
  \node[rtri] (n4) at (c4) {};
  \node[draw, dashed, thick, rounded corners=3pt, inner sep=0.3cm, fit=(n1)(n2)] {};
  \node[draw, dashed, thick, rounded corners=3pt, inner sep=0.3cm, fit=(n4)] {};
  \node at (0.9,0.95) {$P^{(<3)}$};
  \node at (3.6,0.95) {$C^{(3)}$};
  \node at (5.4,0.95) {$Q^{(>3)}$};
\end{tikzpicture}
\end{equation}
and \cref{eq_bond_form} can be represented diagrammatically by
\begin{equation}
\begin{tikzpicture}[baseline=(c1.base)]
  \foreach \i/\x in {1/0,2/1.8,3/3.6,4/5.4} \coordinate (c\i) at (\x,0);
  \coordinate (cS) at (4.5,0);
  \draw[thick] (c1) -- (c4);
  \foreach \i in {1,2,3,4} \draw[thick] (c\i) -- ++(0,-0.8);
  \node[ltri] (n1) at (c1) {};
  \node[ltri] (n2) at (c2) {};
  \node[ltri] (n3) at (c3) {};
  \node[core] (nS) at (cS) {};
  \node[rtri] (n4) at (c4) {};
  \node[draw, dashed, thick, rounded corners=3pt, inner sep=0.3cm, fit=(n1)(n2)(n3)] {};
  \node[draw, dashed, thick, rounded corners=3pt, inner sep=0.3cm, fit=(n4)] {};
  \node at (1.8,0.95) {$P^{(\le 3)}$};
  \node at (4.5,0.95) {$S^{(3)}$};
  \node at (5.4,0.95) {$Q^{(>3)}$};
\end{tikzpicture}
\end{equation}
In this tensor diagram, we use triangles for orthogonal cores.

For $s=2$ ($d=1$), these forms reduce to the familiar three representations
\begin{equation}
\label{eq_Wigner_1D1V_three_forms}
W(x,p)
= \sum_{j=1}^r C_j^{(1)}(x) Q_j^{(2)}(p)
= \sum_{j=1}^r P_j^{(1)}(x) C_j^{(2)}(p)
= \sum_{i=1}^r\sum_{j=1}^r P_i^{(1)}(x) S_{ij}^{(1)} Q_j^{(2)}(p).
\end{equation}
With the TT representation, the memory cost for storing the solution is $\mathcal{O}(d r^2 (N_x+N_p))$, where $r=\max\{r_1,\cdots,r_{s-1}\}$.
For $s=2$, the memory cost further reduces to $\mathcal{O}(r(N_x+N_p+r))$.
When $W$ possesses an intrinsic low-rank structure, the TT rank $r$ can remain relatively small, typically satisfying $r\ll N_x$ and $r\ll N_p$.
The memory cost can therefore be reduced significantly.

\subsection{Low-rank representations of the potential difference}
In addition to the low-rank assumption on the Wigner function, our method assumes that the potential difference $\delta V(\bx,\be)$ admits a suitable low-rank tensor representation.

We mainly consider the canonical polyadic (CP) representation
\begin{equation}
\label{eq_deltaV_CP}
\delta V(\bx,\be) = \sum_{\mu=1}^\varrho
D_\mu^{(1)}(x_1)D_\mu^{(2)}(\eta_1)\cdots
D_\mu^{(s-1)}(x_d)D_\mu^{(s)}(\eta_d),
\end{equation}
where $\varrho$ denotes the number of separable terms in the CP representation.

Alternatively, $\delta V$ can be represented in tensor-train form,
\begin{equation}
\label{eq_deltaV_TT}
\delta V(\bx,\be)
=
\sum_{\nu_1=1}^{\lambda_1}
\sum_{\nu_2=1}^{\lambda_2}
\cdots
\sum_{\nu_{s-1}=1}^{\lambda_{s-1}}
B_{1\nu_1}^{(1)}(x_1)
B_{\nu_1\nu_2}^{(2)}(\eta_1)
\cdots
B_{\nu_{s-2}\nu_{s-1}}^{(s-1)}(x_d)
B_{\nu_{s-1}1}^{(s)}(\eta_d),
\end{equation}
where $(\lambda_1,\ldots,\lambda_{s-1})$ denotes the TT rank of the potential difference.
The TT ranks of $\delta V$ are independent of the prescribed TT ranks $(r_1,\ldots,r_{s-1})$ used to represent the Wigner function $W$.

In the following derivation, we use the CP representation as the primary formulation because of its simpler notation.
The corresponding expressions for a TT representation of $\delta V$ are given whenever the structure of the projected operator is affected.

In the following discussion, we write \cref{eq_wfp} in 
\begin{equation}
\label{eq_wfp_compact}
\partial_t W(t,\bx,\bp) = L_1(W) + L_2(W).
\end{equation}
where $L_1$ and $L_2$ are defined by
\begin{align}
\label{eq_defn_L1}
& L_1(W) = -\frac{\bp}{m}\cdot \nabla_\bx W + 2\gamma \nabla_\bp \cdot (\bp W) + D_{pp} \Delta_\bp W + 2 D_{xp} \nabla_\bx\cdot\nabla_\bp W + D_{xx}\Delta_\bx W, \\
&
\label{eq_defn_L2}
L_2(W) = -\Theta[V] W = \frac{\ii}{\hbar} \fourier^{-1}[\delta V(\bx,\be) \fourier W].
\end{align}
All local operators are collected in $L_1$ and $L_2$ denotes the global Wigner pseudo-differential operator.
The splitting is consistent with the operator splitting techniques for Wigner-type problems \cite{arnold1995operator,arnold1996operator}. 
After the splitting, $L_1$ can be evaluated efficiently in the original phase-space variables, whereas $L_2$ is more naturally evaluated in Fourier space.

We first consider the operator $L_2$, as its treatment is more mathematically involved.
In Fourier space, $L_2$ involves multiplication by the purely imaginary factor $\ii\,\delta V/\hbar$, where $\delta V$ is real-valued.
When the Wigner function is in TT form and the potential difference is in CP form \eqref{eq_deltaV_CP},
$L_2$ admits the following separated representation,
\begin{equation}
\label{eq_L2_CPMPO}
L_2 = \frac{\ii}{\hbar}\sum_{\mu=1}^\varrho
D_\mu^{(1)}(x_1)
\left(\mathcal{F}_{p_1}^{-1}D_\mu^{(2)}(\eta_1)\mathcal{F}_{p_1}\right)
\cdots
D_\mu^{(s-1)}(x_d)
\left(\mathcal{F}_{p_d}^{-1}D_\mu^{(s)}(\eta_d)\mathcal{F}_{p_d}\right)
\end{equation}
where $\mathcal{F}_{p_k}$ and $\mathcal{F}_{p_k}^{-1}$ are the Fourier transform and inverse Fourier transform in the $k$th momentum direction, respectively.

If $\delta V$ is instead represented in a TT form \eqref{eq_deltaV_TT}, the operator can also be written in a similar form:
\begin{equation}
\label{eq_L2_TTMPO}
L_2
=
\frac{\ii}{\hbar}
\sum_{\nu_1=1}^{\lambda_1}
\cdots
\sum_{\nu_{s-1}=1}^{\lambda_{s-1}}
B_{1\nu_1}^{(1)}(x_1)
\left(
\mathcal{F}_{p_1}^{-1}
B_{\nu_1\nu_2}^{(2)}(\eta_1)
\mathcal{F}_{p_1}
\right)
\cdots
B_{\nu_{s-2}\nu_{s-1}}^{(s-1)}(x_d)
\left(
\mathcal{F}_{p_d}^{-1}
B_{\nu_{s-1}1}^{(s)}(\eta_d)
\mathcal{F}_{p_d}
\right).
\end{equation}
Thus, for both representations, each odd factor acts by pointwise multiplication in the corresponding spatial variable $x_q$, whereas each even factor acts through a Fourier transform, pointwise multiplication in $\eta_q$, and an inverse Fourier transform in the corresponding momentum variable $p_q$.

In \cref{eq_L2_CPMPO,eq_L2_TTMPO},
the Fourier transform and inverse Fourier transform with respect to a momentum variable do not increase the TT rank.
The summation and contraction associated with the potential representation, however, may increase the TT ranks substantially.
Therefore, direct application of $L_2$ to $W$, although does not involve full tensor computation, may increase the TT rank of the representation significantly.
The projector-splitting formulation developed below avoids explicitly forming this high-rank intermediate tensor by projecting the action of $L_2$ directly onto the tangent space of the prescribed low-rank TT manifold.
Compared with the treatment of $L_2$, the treatment of $L_1$ is more standard.
Its detailed implementation will be discussed after the spatial and momentum discretizations are introduced in the following sections.

\begin{remark}
Although we focus on CP and TT representations of the potential difference $\delta V$ in this paper, the formulation is not restricted to these two tensor formats.
For example, $\delta V$ may also be represented in Tucker form.
Substituting a Tucker representation of $\delta V$ into \cref{eq_defn_L2} leads to a corresponding tensor-structured representation of the operator $L_2$, analogous to the CP and TT forms in \cref{eq_L2_CPMPO,eq_L2_TTMPO}.
The subsequent treatment can then be developed in a similar manner.
\end{remark}

\section{Projector-splitting method}
\label{sec_projector_splitting_method}
In this section, we apply the tensor-train projector-splitting method \cite{lubich2015time} to the Wigner--Fokker--Planck equation.
A single time step consists of a left-to-right sweep with $s=2d$ center-core steps and $s-1=2d-1$ bond steps, resulting in a total of $4d-1$ local evolution problems.
For $d=1$, the tensor train reduces to a matrix representation, and the resulting scheme coincides with the familiar $K$-, $S$-, and $L$-substep projector-splitting integrator \cite{lubich2014projector}.

We classify the local evolution problems into spatial center-core, momentum center-core, and bond steps.
For each substep, \cref{eq_wfp_compact} is projected onto the corresponding TT subspace.
The projected contribution of $L_2(W)$ is represented by a local matrix operator, while the treatment of $L_1(W)$ is deferred to the next section after the spatial and momentum discretizations are introduced.
We then apply a first-order exponential update to the $L_2$ contribution and a forward Euler update to the $L_1$ contribution.
Finally, a matrix-free implementation of the exponential action is introduced to avoid explicitly forming the local matrix exponential. 

\subsection{Spatial center-core step}
We first consider the spatial center-core steps.
At an odd site $k=2q-1$, we write the solution $W$ in the form of \cref{eq_canonical_form}.
In the projector-splitting method, \cref{eq_wfp_compact} is projected onto the subspace spanned by the left and right basis functions $P_{a_{k-1}}^{(<k)}(\bz_{<k})$ and $Q_{a_k}^{(>k)}(\bz_{>k})$.
We denote the resulting projected equation at $x_q$ by
\begin{equation}
\label{eq_spatial_center_core_difC}
\partial_t C^{(k)}(x_q)
=
F_1^{(k)}[W](x_q) + F_2^{(k)}[W](x_q),
\end{equation}
where $F_1^{(k)}[W]$ and $F_2^{(k)}[W]$ are $r_{k-1}\times r_k$ matrix-valued functions with entries
\begin{align}
\left[F_1^{(k)}[W]\right]_{a_{k-1}a_k}
&=
\left\langle
P_{a_{k-1}}^{(<k)}(\bz_{<k}),
L_1(W),
Q_{a_k}^{(>k)}(\bz_{>k})
\right\rangle_{\bz_{\neq k}}, \\
\left[F_2^{(k)}[W]\right]_{a_{k-1}a_k}
&=
\left\langle
P_{a_{k-1}}^{(<k)}(\bz_{<k}),
L_2(W),
Q_{a_k}^{(>k)}(\bz_{>k})
\right\rangle_{\bz_{\neq k}}.
\end{align}
Here, the inner products are taken over all variables except $z_k=x_q$.
We temporarily retain the first term in its projected form and leave its detailed treatment to the following section.

For the CP representation \cref{eq_deltaV_CP}, the second term can be simplified as
\begin{equation}
\label{eq_spatial_center_core_F2_CP}
F_2^{(k)}[W](x_q)
=
\frac{\ii}{\hbar}\sum_{\mu=1}^{\varrho}
D_\mu^{(k)}(x_q)
\mathsf{L}_\mu^{(k-1)}
C^{(k)}(x_q)
\mathsf{R}_\mu^{(k)}.
\end{equation}
Here, $\mathsf{L}_\mu^{(k-1)}\in\mathbb{C}^{r_{k-1}\times r_{k-1}}$ and $\mathsf{R}_\mu^{(k)}\in\mathbb{C}^{r_k\times r_k}$ denote the left- and right-environment matrices, respectively.
Their entries are given by
\begin{align}
\label{eq_L2_spatial_center_core_left_envorinment}
\left[\mathsf{L}_\mu^{(k-1)}\right]_{ab}
={}&
\Bigg\langle
P_a^{(<k)},
D_\mu^{(1)}
\left(\mathcal{F}_{p_1}^{-1}D_\mu^{(2)}\mathcal{F}_{p_1}\right)
\cdots
D_\mu^{(k-2)}
\left(\mathcal{F}_{p_{q-1}}^{-1}D_\mu^{(k-1)}\mathcal{F}_{p_{q-1}}\right)
P_b^{(<k)}
\Bigg\rangle_{\bz_{<k}},
\end{align}
and
\begin{align}
\label{eq_L2_spatial_center_core_right_envorinment}
\left[\mathsf{R}_\mu^{(k)}\right]_{ab}
={}&
\Bigg\langle
Q_a^{(>k)},
\left(\mathcal{F}_{p_q}^{-1}D_\mu^{(k+1)}\mathcal{F}_{p_q}\right)
D_\mu^{(k+2)}
\cdots
D_\mu^{(s-1)}
\left(\mathcal{F}_{p_d}^{-1}D_\mu^{(s)}\mathcal{F}_{p_d}\right)
Q_b^{(>k)}
\Bigg\rangle_{\bz_{>k}}.
\end{align}
For the left boundary $k=1$, we set $\mathsf{L}_\mu^{(0)}=1$.

For example, when $s=4$ ($d=2$) and $k=3$, the following diagrammatic equation shows the construction of $\mathsf{L}_\mu^{(2)}$ and $\mathsf{R}_\mu^{(3)}$ in $F_2^{(3)}[W]$:
\begin{equation}
F_2^{(3)}[W]= 
\frac{\ii}{\hbar}
\sum_{\mu=1}^{\varrho}
\,
\begin{tikzpicture}[baseline=(mid),
  diam/.style={diamond, draw=black, thick, fill=white, inner sep=0pt, minimum size=0.3cm}]
  \coordinate (mid) at (0,-1.5);
  \foreach \i/\x in {1/0,2/1.8,3/3.6,4/5.4} {
    \coordinate (w\i) at (\x,0);
    \coordinate (m\i) at (\x,-3.2);
  }
  \draw[thick] (w1) -- (w4);
  \draw[thick] (m1) -- (m2) -- (3.3,-3.2);
  \draw[thick] (m4) -- (3.9,-3.2);
  \foreach \i in {1,2,4} \draw[thick] (w\i) -- (m\i);
  \draw[thick] (w3) -- (3.6,-2.4);
  \node[ltri] at (w1) {};
  \node[ltri] at (w2) {};
  \node[core] at (w3) {};
  \node[rtri] at (w4) {};
  \node[ltri] at (m1) {};
  \node[ltri] at (m2) {};
  \node[rtri] at (m4) {};
  \node[fbox] at (1.8,-1.0) {};  \node[right=2pt] at (1.8,-1.0) {$\mathcal F$};
  \node[diam] at (1.8,-1.6) {};  \node[right=2pt] at (1.8,-1.6) {$D^{(2)}_\mu$};
  \node[fbox] at (1.8,-2.2) {};  \node[right=2pt] at (1.8,-2.2) {$\mathcal F^{-1}$};
  \node[fbox] at (5.4,-1.0) {};  \node[right=2pt] at (5.4,-1.0) {$\mathcal F$};
  \node[diam] at (5.4,-1.6) {};  \node[right=2pt] at (5.4,-1.6) {$D^{(4)}_\mu$};
  \node[fbox] at (5.4,-2.2) {};  \node[right=2pt] at (5.4,-2.2) {$\mathcal F^{-1}$};
  \node[diam] at (0,-1.6)   {};  \node[right=2pt] at (0,-1.6)   {$D^{(1)}_\mu$};
  \node[diam] at (3.6,-1.6) {};  \node[right=2pt] at (3.6,-1.6) {$D^{(3)}_\mu$};
  \draw[dashed, thick, rounded corners=3pt] (-0.5,0.6) rectangle (3.0,-3.8);
  \draw[dashed, thick, rounded corners=3pt] (4.8,0.6) rectangle (6.7,-3.8);
  \node at (1.25,0.95) {$\mathsf{L}^{(2)}_\mu$};
  \node at (3.6,0.95)  {$C^{(3)}$};
  \node at (5.75,0.95) {$\mathsf{R}^{(3)}_\mu$};
\end{tikzpicture}\quad.
\end{equation}

The environment matrices $\mathsf{L}_\mu^{(k-1)}$ and $\mathsf{R}_\mu^{(k)}$ therefore collect the partial contractions of the $\mu$th separable term of the potential operator with all TT modes to the left and right of the current center core, respectively.
After these contractions are computed, the dependence on the remaining physical variable $z_k=x_q$ is contained entirely in the local factor $D_\mu^{(k)}(x_q)$ and the center core $C^{(k)}(x_q)$.

Using \cref{eq_spatial_center_core_F2_CP} and the identity
$\vec{ACB}=(B^\top\otimes A)\vec{C}$,
\cref{eq_spatial_center_core_difC} can be written as
\begin{equation}
\label{eq_spatial_center_core_vectorized}
\partial_t \vec{C^{(k)}(x_q)}
=
\vec{F_1^{(k)}[W](x_q)}
+
M_k(x_q)\vec{C^{(k)}(x_q)},
\end{equation}
where
\begin{equation}
\label{eq_spatial_M_CP}
M_k(x_q)
=
\frac{\ii}{\hbar}
\sum_{\mu=1}^{\varrho}
D_\mu^{(k)}(x_q)
\left(\mathsf{R}_\mu^{(k)}\right)^\top
\otimes
\mathsf{L}_\mu^{(k-1)}.
\end{equation}

When the potential difference is represented in TT form \cref{eq_deltaV_TT}, the same projection procedure gives
\begin{equation}
\label{eq_spatial_center_core_F2_TT}
F_2^{(k)}[W](x_q)
=
\frac{\ii}{\hbar}
\sum_{\nu_{k-1}=1}^{\lambda_{k-1}}
\sum_{\nu_k=1}^{\lambda_k}
B_{\nu_{k-1}\nu_k}^{(k)}(x_q)
\mathsf{L}_{\nu_{k-1}}^{(k-1)}
C^{(k)}(x_q)
\mathsf{R}_{\nu_k}^{(k)}.
\end{equation}
Here,
$\mathsf{L}_{\nu_{k-1}}^{(k-1)}
\in\mathbb{C}^{r_{k-1}\times r_{k-1}}$
and
$\mathsf{R}_{\nu_k}^{(k)}
\in\mathbb{C}^{r_k\times r_k}$
denote the left- and right-environment matrices for the TT representation of the potential difference, respectively.
Their entries are given by
\begin{align}
\label{eq_L2_spatial_center_core_left_environment_TT}
\left[\mathsf{L}_{\nu_{k-1}}^{(k-1)}\right]_{ab}
={}&
\sum_{\nu_1=1}^{\lambda_1}
\cdots
\sum_{\nu_{k-2}=1}^{\lambda_{k-2}}
\Bigg\langle
P_a^{(<k)},
B_{1\nu_1}^{(1)}
\left(
\mathcal{F}_{p_1}^{-1}
B_{\nu_1\nu_2}^{(2)}
\mathcal{F}_{p_1}
\right)
\cdots
B_{\nu_{k-3}\nu_{k-2}}^{(k-2)}
\left(
\mathcal{F}_{p_{q-1}}^{-1}
B_{\nu_{k-2}\nu_{k-1}}^{(k-1)}
\mathcal{F}_{p_{q-1}}
\right)
P_b^{(<k)}
\Bigg\rangle_{\bz_{<k}},
\end{align}
and
\begin{align}
\label{eq_L2_spatial_center_core_right_environment_TT}
\left[\mathsf{R}_{\nu_k}^{(k)}\right]_{ab}
={}&
\sum_{\nu_{k+1}=1}^{\lambda_{k+1}}
\cdots
\sum_{\nu_{s-1}=1}^{\lambda_{s-1}}
\Bigg\langle
Q_a^{(>k)},
\left(
\mathcal{F}_{p_q}^{-1}
B_{\nu_k\nu_{k+1}}^{(k+1)}
\mathcal{F}_{p_q}
\right)
B_{\nu_{k+1}\nu_{k+2}}^{(k+2)}
\cdots
B_{\nu_{s-2}\nu_{s-1}}^{(s-1)}
\left(
\mathcal{F}_{p_d}^{-1}
B_{\nu_{s-1}1}^{(s)}
\mathcal{F}_{p_d}
\right)
Q_b^{(>k)}
\Bigg\rangle_{\bz_{>k}}.
\end{align}
For the left boundary $k=1$, we set $\mathsf{L}_{1}^{(0)}=1$.
Empty sums and products in the boundary cases are understood in the usual sense.

Similar to the CP case, the environment matrices
$\mathsf{L}_{\nu_{k-1}}^{(k-1)}$
and
$\mathsf{R}_{\nu_k}^{(k)}$
collect the partial contractions of the potential TT cores with all solution TT modes to the left and right of the current center core, respectively.
The potential TT indices $\nu_{k-1}$ and $\nu_k$ remain uncontracted and connect the two environments to the current potential core
$B_{\nu_{k-1}\nu_k}^{(k)}(x_q)$.
After the left and right contractions are computed, the dependence on the remaining physical variable $z_k=x_q$ is contained entirely in
$B_{\nu_{k-1}\nu_k}^{(k)}(x_q)$
and the center core $C^{(k)}(x_q)$.

The vectorized equation retains the form \cref{eq_spatial_center_core_vectorized}, with the local matrix now given by
\begin{equation}
\label{eq_spatial_M_TT}
M_k(x_q)
=
\frac{\ii}{\hbar}
\sum_{\nu_{k-1}=1}^{\lambda_{k-1}}
\sum_{\nu_k=1}^{\lambda_k}
B_{\nu_{k-1}\nu_k}^{(k)}(x_q)
\left(\mathsf{R}_{\nu_k}^{(k)}\right)^\top
\otimes
\mathsf{L}_{\nu_{k-1}}^{(k-1)}.
\end{equation}
Thus, the CP and TT representations lead to the same local evolution equation.
The difference lies in the structure of the potential environments: in the CP representation, the same separable index $\mu$ is shared by the left and right contractions, whereas in the TT representation, the two contractions are connected to the current potential core through the neighboring TT indices $\nu_{k-1}$ and $\nu_k$.

\subsection{Momentum center-core step}

At an even site $k=2q$, we again write the solution $W$ in the form of \cref{eq_canonical_form}.
Following the same projection procedure as in the spatial center-core step, we obtain
\begin{equation}
\label{eq_momentum_center_core_difC}
\partial_t C^{(k)}(p_q)
=
F_1^{(k)}[W](p_q) + F_2^{(k)}[W](p_q),
\end{equation}
where $F_1^{(k)}[W]$ and $F_2^{(k)}[W]$ are defined by the same partial projections as in the spatial center-core step, with the inner products now taken over all variables except $z_k=p_q$.
As before, we temporarily retain $F_1^{(k)}[W]$ in its projected form and defer its detailed treatment to the following section.

For the CP representation, the $k$th factor of the operator acts in Fourier space, giving
\begin{equation}
\label{eq_momentum_center_core_F2_CP}
F_2^{(k)}[W](p_q)
=
\frac{\ii}{\hbar}
\sum_{\mu=1}^{\varrho}
\mathsf{L}_\mu^{(k-1)}
\mathcal{F}_{p_q}^{-1}
\left[
D_\mu^{(k)}(\eta_q)
\mathcal{F}_{p_q}C^{(k)}(p_q)
\right]
\mathsf{R}_\mu^{(k)}.
\end{equation}
Here,
$\mathsf{L}_\mu^{(k-1)}
\in\mathbb{C}^{r_{k-1}\times r_{k-1}}$
and
$\mathsf{R}_\mu^{(k)}
\in\mathbb{C}^{r_k\times r_k}$
denote the left- and right-environment matrices, respectively.
Their entries are
\begin{align}
\label{eq_L2_momentum_center_core_left_environment}
\left[\mathsf{L}_\mu^{(k-1)}\right]_{ab}
={}&
\Bigg\langle
P_a^{(<k)},
D_\mu^{(1)}
\left(
\mathcal{F}_{p_1}^{-1}
D_\mu^{(2)}
\mathcal{F}_{p_1}
\right)
\cdots
\left(
\mathcal{F}_{p_{q-1}}^{-1}
D_\mu^{(k-2)}
\mathcal{F}_{p_{q-1}}
\right)
D_\mu^{(k-1)}
P_b^{(<k)}
\Bigg\rangle_{\bz_{<k}},
\end{align}
and
\begin{align}
\label{eq_L2_momentum_center_core_right_environment}
\left[\mathsf{R}_\mu^{(k)}\right]_{ab}
={}&
\Bigg\langle
Q_a^{(>k)},
D_\mu^{(k+1)}
\left(
\mathcal{F}_{p_{q+1}}^{-1}
D_\mu^{(k+2)}
\mathcal{F}_{p_{q+1}}
\right)
\cdots
D_\mu^{(s-1)}
\left(
\mathcal{F}_{p_d}^{-1}
D_\mu^{(s)}
\mathcal{F}_{p_d}
\right)
Q_b^{(>k)}
\Bigg\rangle_{\bz_{>k}}.
\end{align}
For the right boundary $k=s$, we set $\mathsf{R}_\mu^{(s)}=1$.

For example, when $s=4$ ($d=2$) and $k=2$, the following diagrammatic equation shows the construction of $\mathsf{L}_\mu^{(1)}$ and $\mathsf{R}_\mu^{(2)}$ in $F_2^{(2)}[W]$:
\begin{equation}
F_2^{(2)}[W] = \frac{\ii}{\hbar} \sum_{\mu=1}^{\varrho} \,
\begin{tikzpicture}[baseline=(mid),
  diam/.style={diamond, draw=black, thick, fill=white, inner sep=0pt, minimum size=0.3cm}]
  \coordinate (mid) at (0,-1.5);
  \foreach \i/\x in {1/0,2/1.8,3/3.6,4/5.4} {
    \coordinate (w\i) at (\x,0);
    \coordinate (m\i) at (\x,-3.2);
  }
  \draw[thick] (w1) -- (w4);
  \draw[thick] (m1) -- (1.5,-3.2);
  \draw[thick] (m3) -- (m4);
  \draw[thick] (m3) -- (2.1,-3.2);
  \foreach \i in {1,3,4} \draw[thick] (w\i) -- (m\i);
  \draw[thick] (w2) -- (1.8,-2.8);
  \node[ltri] at (w1) {};
  \node[core] at (w2) {};
  \node[rtri] at (w3) {};
  \node[rtri] at (w4) {};
  \node[ltri] at (m1) {};
  \node[rtri] at (m3) {};
  \node[rtri] at (m4) {};
  \node[diam] at (0,-1.6)   {};  \node[right=2pt] at (0,-1.6)   {$D^{(1)}_\mu$};
  \node[fbox] at (1.8,-1.0) {};  \node[right=2pt] at (1.8,-1.0) {$\mathcal F$};
  \node[diam] at (1.8,-1.6) {};  \node[right=2pt] at (1.8,-1.6) {$D^{(2)}_\mu$};
  \node[fbox] at (1.8,-2.2) {};  \node[right=2pt] at (1.8,-2.2) {$\mathcal F^{-1}$};
  \node[diam] at (3.6,-1.6) {};  \node[right=2pt] at (3.6,-1.6) {$D^{(3)}_\mu$};
  \node[fbox] at (5.4,-1.0) {};  \node[right=2pt] at (5.4,-1.0) {$\mathcal F$};
  \node[diam] at (5.4,-1.6) {};  \node[right=2pt] at (5.4,-1.6) {$D^{(4)}_\mu$};
  \node[fbox] at (5.4,-2.2) {};  \node[right=2pt] at (5.4,-2.2) {$\mathcal F^{-1}$};
  \draw[dashed, thick, rounded corners=3pt] (-0.5,0.6) rectangle (1.2,-3.8);
  \draw[dashed, thick, rounded corners=3pt] (3.0,0.6) rectangle (6.7,-3.8);
  \node at (0.35,0.95) {$\mathsf{L}^{(1)}_\mu$};
  \node at (1.8,0.95)  {$C^{(2)}$};
  \node at (4.85,0.95) {$\mathsf{R}^{(2)}_\mu$};
\end{tikzpicture}\quad.
\end{equation}

The environment matrices
$\mathsf{L}_\mu^{(k-1)}$
and
$\mathsf{R}_\mu^{(k)}$
collect the partial contractions of the $\mu$th separable term of the potential operator with all TT modes to the left and right of the current momentum center core, respectively.
After these contractions are computed, the dependence on the remaining physical variable is contained in the local factor
$D_\mu^{(k)}(\eta_q)$
and the center core in Fourier space.

To expose the local structure of the evolution, we define
\begin{equation}
\widehat{C}^{(k)}(\eta_q)
=
\mathcal{F}_{p_q}C^{(k)}(p_q).
\end{equation}
Applying $\mathcal{F}_{p_q}$ to \cref{eq_momentum_center_core_difC} and using \cref{eq_momentum_center_core_F2_CP} gives
\begin{equation}
\label{eq_momentum_center_core_vectorized}
\partial_t \vec{\widehat{C}^{(k)}(\eta_q)}
=
\vec{\mathcal{F}_{p_q}\left(F_1^{(k)}[W]\right)(\eta_q)}
+
M_k(\eta_q)
\vec{\widehat{C}^{(k)}(\eta_q)},
\end{equation}
where
\begin{equation}
\label{eq_momentum_M_CP}
M_k(\eta_q)
=
\frac{\ii}{\hbar}
\sum_{\mu=1}^{\varrho}
D_\mu^{(k)}(\eta_q)
\left(\mathsf{R}_\mu^{(k)}\right)^\top
\otimes
\mathsf{L}_\mu^{(k-1)}.
\end{equation}

When $\delta V$ is represented in TT form, the projected $L_2$ contribution become
\begin{equation}
\label{eq_momentum_center_core_F2_TT}
F_2^{(k)}[W](p_q)
=
\frac{\ii}{\hbar}
\sum_{\nu_{k-1}=1}^{\lambda_{k-1}}
\sum_{\nu_k=1}^{\lambda_k}
\mathsf{L}_{\nu_{k-1}}^{(k-1)}
\mathcal{F}_{p_q}^{-1}
\left[
B_{\nu_{k-1}\nu_k}^{(k)}(\eta_q)
\mathcal{F}_{p_q}C^{(k)}(p_q)
\right]
\mathsf{R}_{\nu_k}^{(k)}.
\end{equation}
Here,
$\mathsf{L}_{\nu_{k-1}}^{(k-1)}
\in\mathbb{C}^{r_{k-1}\times r_{k-1}}$
and
$\mathsf{R}_{\nu_k}^{(k)}
\in\mathbb{C}^{r_k\times r_k}$
denote the corresponding left- and right-environment matrices for the TT representation of the potential difference.
Their entries are
\begin{align}
\label{eq_L2_momentum_center_core_left_environment_TT}
\left[\mathsf{L}_{\nu_{k-1}}^{(k-1)}\right]_{ab}
={}&
\sum_{\nu_1=1}^{\lambda_1}
\cdots
\sum_{\nu_{k-2}=1}^{\lambda_{k-2}}
\Bigg\langle
P_a^{(<k)},
B_{1\nu_1}^{(1)}
\left(
\mathcal{F}_{p_1}^{-1}
B_{\nu_1\nu_2}^{(2)}
\mathcal{F}_{p_1}
\right)
\cdots
\left(
\mathcal{F}_{p_{q-1}}^{-1}
B_{\nu_{k-3}\nu_{k-2}}^{(k-2)}
\mathcal{F}_{p_{q-1}}
\right)
B_{\nu_{k-2}\nu_{k-1}}^{(k-1)}
P_b^{(<k)}
\Bigg\rangle_{\bz_{<k}},
\end{align}
and
\begin{align}
\label{eq_L2_momentum_center_core_right_environment_TT}
\left[\mathsf{R}_{\nu_k}^{(k)}\right]_{ab}
={}&
\sum_{\nu_{k+1}=1}^{\lambda_{k+1}}
\cdots
\sum_{\nu_{s-1}=1}^{\lambda_{s-1}}
\Bigg\langle
Q_a^{(>k)},
B_{\nu_k\nu_{k+1}}^{(k+1)}
\left(
\mathcal{F}_{p_{q+1}}^{-1}
B_{\nu_{k+1}\nu_{k+2}}^{(k+2)}
\mathcal{F}_{p_{q+1}}
\right)
\cdots
B_{\nu_{s-2}\nu_{s-1}}^{(s-1)}
\left(
\mathcal{F}_{p_d}^{-1}
B_{\nu_{s-1}1}^{(s)}
\mathcal{F}_{p_d}
\right)
Q_b^{(>k)}
\Bigg\rangle_{\bz_{>k}}.
\end{align}
For the right boundary $k=s$, we set $\mathsf{R}_{1}^{(s)}=1$.
Empty sums and products in the boundary cases are understood in the usual sense.

Similar to the CP case, the environment matrices
$\mathsf{L}_{\nu_{k-1}}^{(k-1)}$
and
$\mathsf{R}_{\nu_k}^{(k)}$
collect the partial contractions of the potential TT cores with all solution TT modes to the left and right of the current momentum center core, respectively.
The potential TT indices $\nu_{k-1}$ and $\nu_k$ remain uncontracted and connect the two environments to the current potential core
$B_{\nu_{k-1}\nu_k}^{(k)}(\eta_q)$.
After the left and right contractions are computed, the dependence on the remaining Fourier variable $\eta_q$ is contained entirely in
$B_{\nu_{k-1}\nu_k}^{(k)}(\eta_q)$
and
$\widehat{C}^{(k)}(\eta_q)$.

The Fourier-transformed center-core equation retains the form \cref{eq_momentum_center_core_vectorized}, with
\begin{equation}
\label{eq_momentum_M_TT}
M_k(\eta_q)
=
\frac{\ii}{\hbar}
\sum_{\nu_{k-1}=1}^{\lambda_{k-1}}
\sum_{\nu_k=1}^{\lambda_k}
B_{\nu_{k-1}\nu_k}^{(k)}(\eta_q)
\left(\mathsf{R}_{\nu_k}^{(k)}\right)^\top
\otimes
\mathsf{L}_{\nu_{k-1}}^{(k-1)}.
\end{equation}
Thus, the momentum center-core evolution has the same local matrix structure as the spatial center-core evolution, except that the local evolution associated with $L_2$ is carried out in Fourier space.

\subsection{Bond step}
After each center-core step at site $k$, except for the last site $k=s$, the projector-splitting method performs a bond step between sites $k$ and $k+1$.
We write the solution $W$ in the bond form \cref{eq_bond_form}.
An important distinction from the center-core steps is that the corresponding subprojector in the projector-splitting decomposition carries a negative sign \cite{lubich2014projector,lubich2015time}.
Consequently, the local evolution equation for the bond matrix $S^{(k)}$ is
\begin{equation}
\label{eq_bond_difS}
\partial_t S^{(k)}
=
-G_1^{(k)}[W]
-G_2^{(k)}[W],
\end{equation}
where $G_1^{(k)}[W]$ and $G_2^{(k)}[W]$ are $r_k\times r_k$ matrices obtained by projecting $L_1(W)$ and $L_2(W)$, respectively, onto the left and right basis functions in \cref{eq_bond_form}.
Their entries are
\begin{align}
\left[G_1^{(k)}[W]\right]_{ab}
&=
\left\langle
P_a^{(\leqslant k)}(\bz_{\leqslant k}),
L_1(W),
Q_b^{(>k)}(\bz_{>k})
\right\rangle_{\bz}, \\
\left[G_2^{(k)}[W]\right]_{ab}
&=
\left\langle
P_a^{(\leqslant k)}(\bz_{\leqslant k}),
L_2(W),
Q_b^{(>k)}(\bz_{>k})
\right\rangle_{\bz}.
\end{align}
Here, the inner products are taken over all phase-space variables, since the bond matrix $S^{(k)}$ does not depend on any physical variable.

As in the center-core steps, we temporarily retain $G_1^{(k)}[W]$ in its projected form and defer its detailed treatment to the following section.

For the CP representation, contracting the separated potential operator with the left and right basis functions gives
\begin{equation}
\label{eq_bond_G2_CP}
G_2^{(k)}[W]
=
\frac{\ii}{\hbar}
\sum_{\mu=1}^{\varrho}
\mathsf{L}_\mu^{(k)}
S^{(k)}
\mathsf{R}_\mu^{(k)}.
\end{equation}
Here,
$\mathsf{L}_\mu^{(k)}
\in\mathbb{C}^{r_k\times r_k}$
and
$\mathsf{R}_\mu^{(k)}
\in\mathbb{C}^{r_k\times r_k}$
denote the left- and right-environment matrices associated with the bond between sites $k$ and $k+1$, respectively.
Similar to the center-core steps, these matrices collect the partial contractions of the $\mu$th separable term of the potential operator with all TT modes on the two sides of the bond.

For example, when $s=4$ ($d=2$) and $k=3$, the following diagrammatic equation of $G_2^{(3)}[W]$ shows the construction of $\mathsf{L}_\mu^{(3)}$ and $\mathsf{R}_\mu^{(3)}$:
\begin{equation}
G_2^{(3)}[W] = \frac{\ii}{\hbar} \sum_{\mu=1}^\varrho \,
\begin{tikzpicture}[baseline=(mid),
  diam/.style={diamond, draw=black, thick, fill=white, inner sep=0pt, minimum size=0.3cm}]
  \coordinate (mid) at (0,-1.5);
  \foreach \i/\x in {1/0,2/2.1,3/4.2,4/6.3} {
    \coordinate (w\i) at (\x,0);
    \coordinate (m\i) at (\x,-3.2);
  }
  \coordinate (wS) at (5.25,0);
  \draw[thick] (w1) -- (w4);
  \draw[thick] (m1) -- (m3) -- (4.9,-3.2);
  \draw[thick] (m4) -- (5.6,-3.2);
  \foreach \i in {1,2,3,4} \draw[thick] (w\i) -- (m\i);
  \node[ltri] at (w1) {};
  \node[ltri] at (w2) {};
  \node[ltri] at (w3) {};
  \node[core] at (wS) {};
  \node[rtri] at (w4) {};
  \node[ltri] at (m1) {};
  \node[ltri] at (m2) {};
  \node[ltri] at (m3) {};
  \node[rtri] at (m4) {};
  \node[diam] at (0,-1.6)   {};  \node[right=2pt] at (0,-1.6)   {$D^{(1)}_\mu$};
  \node[fbox] at (2.1,-1.0) {};  \node[right=2pt] at (2.1,-1.0) {$\mathcal F$};
  \node[diam] at (2.1,-1.6) {};  \node[right=2pt] at (2.1,-1.6) {$D^{(2)}_\mu$};
  \node[fbox] at (2.1,-2.2) {};  \node[right=2pt] at (2.1,-2.2) {$\mathcal F^{-1}$};
  \node[diam] at (4.2,-1.6) {};  \node[left=2pt]  at (4.2,-1.6) {$D^{(3)}_\mu$};
  \node[fbox] at (6.3,-1.0) {};  \node[right=2pt] at (6.3,-1.0) {$\mathcal F$};
  \node[diam] at (6.3,-1.6) {};  \node[right=2pt] at (6.3,-1.6) {$D^{(4)}_\mu$};
  \node[fbox] at (6.3,-2.2) {};  \node[right=2pt] at (6.3,-2.2) {$\mathcal F^{-1}$};
  \draw[dashed, thick, rounded corners=3pt] (-0.5,0.6) rectangle (4.65,-3.8);
  \draw[dashed, thick, rounded corners=3pt] (5.85,0.6) rectangle (7.5,-3.8);
  \node at (2.075,0.95) {$\mathsf{L}^{(3)}_\mu$};
  \node at (5.25,0.95)  {$S^{(3)}$};
  \node at (6.675,0.95) {$\mathsf{R}^{(3)}_\mu$};
\end{tikzpicture}\quad.
\end{equation}

Using \cref{eq_bond_G2_CP}, \cref{eq_bond_difS} can be written as
\begin{equation}
\label{eq_bond_vectorized}
\partial_t \vec{S^{(k)}}
=
-\vec{G_1^{(k)}[W]}
-
M_k^{(S)}\vec{S^{(k)}},
\end{equation}
where
\begin{equation}
\label{eq_bond_M_CP}
M_k^{(S)}
=
\frac{\ii}{\hbar}
\sum_{\mu=1}^{\varrho}
\left(\mathsf{R}_\mu^{(k)}\right)^\top
\otimes
\mathsf{L}_\mu^{(k)}.
\end{equation}

When the potential difference is represented in TT form, the potential bond index $\nu_k$ remains uncontracted across the cut between sites $k$ and $k+1$.
The projected contribution becomes
\begin{equation}
\label{eq_bond_G2_TT}
G_2^{(k)}[W]
=
\frac{\ii}{\hbar}
\sum_{\nu_k=1}^{\lambda_k}
\mathsf{L}_{\nu_k}^{(k)}
S^{(k)}
\mathsf{R}_{\nu_k}^{(k)}.
\end{equation}
Here,
$\mathsf{L}_{\nu_k}^{(k)}
\in\mathbb{C}^{r_k\times r_k}$
and
$\mathsf{R}_{\nu_k}^{(k)}
\in\mathbb{C}^{r_k\times r_k}$
denote the corresponding left- and right-environment matrices for the TT representation of the potential difference.
They are obtained by contracting all potential TT cores and solution TT modes to the left and right of the $k$th bond, while keeping the potential bond index $\nu_k$ uncontracted.
Thus, $\nu_k$ connects the two environments across the current bond in the same way that the index $\mu$ identifies the corresponding separable term in the CP representation.
The vectorized bond equation retains the form \cref{eq_bond_vectorized}:
\begin{equation}
\label{eq_bond_M_TT}
    M_{k}^{(S)} = \frac{\ii}{\hbar} \sum_{\nu_k=1}^{\lambda_k}\left(\mathsf{R}_{\nu_k}^{(k)}\right)^\top \otimes \mathsf{L}_{\nu_k}^{(k)}.
\end{equation}

\subsection{First-order time discretization}

We next introduce the time discretization for the local projector-splitting subproblems.
The two contributions $L_1$ and $L_2$ have substantially different numerical structures, and we therefore treat them differently in time.

In Fourier space, the Wigner pseudo-differential operator $L_2$ acts through multiplication by the purely imaginary factor $\ii\,\delta V/\hbar$.
Consequently, applying the forward Euler method directly to this contribution leads to numerical instability, since the corresponding linear operator has eigenvalues along the imaginary axis.
Higher-order explicit Runge--Kutta methods can provide larger stability regions, but they require multiple evaluations of the projected $L_2$ operator within each substep.
Moreover, since the projector-splitting integrator considered here is first order in time, using a higher-order Runge--Kutta method for the individual subproblems would not increase the overall temporal order of the scheme.

On the other hand, after projection onto a center core or bond, the $L_2$ contribution is linear and represented by the local matrices derived above.
Its evolution can therefore be integrated directly through the matrix exponential.
For the $L_1$ contribution, we use the forward Euler method, which is consistent with the first-order accuracy of the projector-splitting scheme.
This leads to a first-order exponential-type integrator in which the projected $L_2$ contribution is treated by an exponential update and the projected $L_1$ contribution by an explicit Euler increment.

The following update formulas apply to both representations of the potential difference.
For the CP representation, $M_k(x_q)$, $M_k(\eta_q)$, and $M_k^{(S)}$ are defined by
\cref{eq_spatial_M_CP,eq_momentum_M_CP,eq_bond_M_CP}, respectively.
For the TT representation, they are replaced by
\cref{eq_spatial_M_TT,eq_momentum_M_TT,eq_bond_M_TT}.
The form of the local time integration remains unchanged.

For a spatial center-core step, we use the following update at each spatial point $x_q$:
\begin{equation}
\label{eq_spatial_center_core_update}
C^{(k),\mathrm{out}}(x_q)
=
\vecinv{
\e^{M_k(x_q)\dt}
\vec{C^{(k),\mathrm{in}}(x_q)}
}
+
\dt\left(F_1^{(k)}[W^{\mathrm{in}}]\right)(x_q).
\end{equation}

For a momentum center-core step, we first evolve the core in Fourier space and then transform it back to momentum space:
\begin{equation}
\label{eq_momentum_center_core_update}
C^{(k),\mathrm{out}}(p_q)
=
\mathcal{F}_{p_q}^{-1}
\left(
\vecinv{
\e^{M_k(\eta_q)\dt}
\vec{\widehat{C}^{(k),\mathrm{in}}(\eta_q)}
}
\right)(p_q)
+
\dt\left(F_1^{(k)}[W^{\mathrm{in}}]\right)(p_q).
\end{equation}

For a bond step, the corresponding update is
\begin{equation}
\label{eq_bond_core_update}
S^{(k),\mathrm{out}}
=
\vecinv{
\e^{-M_k^{(S)}\dt}
\vec{S^{(k),\mathrm{in}}}
}
-
\dt G_1^{(k)}[W^{\mathrm{in}}].
\end{equation}

The preceding discussion describes the individual local subproblems in the projector-splitting method.
To summarize the overall procedure, a single time step of the first-order projector-splitting method consists of a left-to-right sweep in the following order:
\begin{equation}
C^{(1)} \to S^{(1)} \to C^{(2)} \to S^{(2)}
\to \cdots \to S^{(s-1)} \to C^{(s)}.
\end{equation}
Thus, each time step contains a total of $2s-1=4d-1$ substeps: $d$ spatial center-core steps, $d$ momentum center-core steps, and $2d-1$ bond steps.

\subsection{Matrix-free implementation}
\label{sec_matrix_free_implementation}

\Cref{eq_spatial_center_core_update,eq_momentum_center_core_update,eq_bond_core_update} provide compact expressions for the three types of projector-splitting substeps.
A direct implementation, however, would require explicitly assembling the local matrices $M_k(x_q)$, $M_k(\eta_q)$, or $M_k^{(S)}$ and computing their matrix exponentials.
For a uniform solution TT rank $r$, the local matrices associated with a generic interior subproblem have dimensions of order $r^2\times r^2$.
Treating them as general dense matrices therefore requires $\mathcal{O}(r^4)$ storage, while a dense matrix exponential generally incurs a computational cost of $\mathcal{O}(r^6)$.
Such an implementation can become prohibitively expensive as the prescribed TT rank increases.

In practice, we avoid explicitly forming both the local matrices and their exponentials.
We use the spatial center-core step as an example, while the momentum center-core and bond steps are treated analogously.
In \cref{eq_spatial_center_core_update}, the matrix $\e^{M_k(x_q)\dt}$ itself is not required; only its action
$
\e^{M_k(x_q)\dt}\vec{C^{(k),\mathrm{in}}(x_q)}
$
needs to be evaluated.
For simplicity, we denote this operation by $\e^{M\dt}\vec{C}$ in the following discussion.

Our implementation uses the exponential-action algorithm in \cite{almohy2011computing}, implemented by 
\texttt{expm\_multiply} in
\texttt{scipy.sparse.linalg}.
For the CP representation of the potential difference, the action of $M$ can also be evaluated without explicitly assembling the corresponding $r^2\times r^2$ matrix.
For example, according to \cref{eq_spatial_center_core_F2_CP},
\begin{equation}
M_k(x_q)\vec{C}
=
\frac{\ii}{\hbar}
\vec{
\sum_{\mu=1}^{\varrho}
D_\mu^{(k)}(x_q)
\mathsf{L}_\mu^{(k-1)}
C
\mathsf{R}_\mu^{(k)}
}.
\end{equation}
Each term involves matrix products of size at most $r\times r$.
Hence, one application of the local operator costs $\mathcal{O}(\varrho r^3)$ operations.
If the exponential-action algorithm requires a total of $n_{\mathrm{mv}}$ such operator applications, the corresponding computational cost is
$\mathcal{O}(n_{\mathrm{mv}}\varrho r^3)$.
The number $n_{\mathrm{mv}}$ depends on the scaling parameter, the Taylor degree, the time step, and the norm of the local operator.
When $n_{\mathrm{mv}}$ remains moderate, the matrix-free formulation reduces the rank dependence substantially compared with the $\mathcal{O}(r^6)$ cost of computing a dense matrix exponential and simultaneously avoids the $\mathcal{O}(r^4)$ storage required for the local matrix.

The same matrix-free strategy applies when $\delta V$ is represented in TT form.
For a center-core step, \cref{eq_spatial_center_core_F2_TT,eq_momentum_center_core_F2_TT} show that one application of the local operator involves summation over the two neighboring potential TT indices $\nu_{k-1}$ and $\nu_k$.
For a uniform solution rank $r$, the corresponding cost is
$\mathcal{O}(\lambda_{k-1}\lambda_k r^3)$.
Defining
$\lambda=\max_j\lambda_j$,
this cost is bounded by $\mathcal{O}(\lambda^2r^3)$.
For a bond step, only the potential bond index $\nu_k$ remains, and one operator application costs
$\mathcal{O}(\lambda_k r^3)$, bounded by $\mathcal{O}(\lambda r^3)$.

Therefore, the same exponential-action algorithm can be used for both CP and TT representations of the potential difference without explicitly forming the local matrices.
The difference lies only in the cost of each operator application: it depends linearly on the number of CP terms $\varrho$ in the CP representation and on the neighboring potential TT ranks in the TT representation.

\section{Spatial and momentum discretization}
\label{sec_discretization}

In this section, we introduce the discretization of the spatial, momentum, and Fourier variables and describe how the operator $L_1$ is represented in the discrete tensor-train setting.

We truncate the spatial and momentum domains to
$\Omega_\bx=[x_{\min},x_{\max}]^d$
and
$\Omega_\bp=[p_{\min},p_{\max}]^d$,
respectively.
For simplicity, we use the same interval in each spatial direction and the same interval in each momentum direction, although this is not required by the method.
In practice, the truncation bounds are chosen sufficiently large so that the computational domain captures essentially all of the Wigner function over the time interval of interest.
In other words, $W(t,\bx,\bp)$ is assumed to be negligible outside $\Omega_\bx\times\Omega_\bp$.

Each spatial direction is divided uniformly into $N_x$ cells and each momentum direction into $N_p$ cells, with mesh sizes
\begin{equation}
\Delta x=\frac{x_{\max}-x_{\min}}{N_x},\qquad
\Delta p=\frac{p_{\max}-p_{\min}}{N_p}.
\end{equation}
We use cell-centered nodes,
\begin{equation}
x_{q,i}=x_{\min}+\left(i-\frac12\right)\Delta x,\qquad
p_{q,j}=p_{\min}+\left(j-\frac12\right)\Delta p,
\end{equation}
for $q=1,\ldots,d$, $i=1,\ldots,N_x$, and $j=1,\ldots,N_p$.
Hence, for fixed $t$, the discrete Wigner function is a $2d$-way tensor of size
$N_x\times N_p\times\cdots\times N_x\times N_p$,
with the same variable ordering as that used in the TT representation.
In our method, a prescribed TT rank
$(r_1,r_2,\cdots,r_{2d-1})$
is fixed, and the full tensor is approximated by a tensor train, substantially reducing the memory requirement.

We consider the discretization of $L_1$.
The operator is naturally decomposed into contributions associated with each pair $(x_q,p_q)$,
\begin{equation}
L_1(W)=\sum_{q=1}^d L_{1,q}(W),
\end{equation}
where
\begin{equation}
L_{1,q}(W)
=
-\frac{p_q}{m}\partial_{x_q}W
+2\gamma\partial_{p_q}(p_qW)
+D_{pp}\partial_{p_qp_q}W
+2D_{xp}\partial_{x_qp_q}W
+D_{xx}\partial_{x_qx_q}W.
\end{equation}
We discretize these terms in conservative flux form.
This choice is convenient both numerically and algebraically. The transport, diffusion, and drift contributions remain local in the corresponding physical variables, and their action can therefore be incorporated directly into the associated TT cores.

Suppressing all indices other than those associated with $(x_q,p_q)$, the discrete operator is written as
\begin{equation}
(L_{1,q}W)_{i,j}
=
-\frac{\mathcal{T}_{i+\frac12,j}-\mathcal{T}_{i-\frac12,j}}{\Delta x}
+\frac{\mathcal{X}_{i+\frac12,j}-\mathcal{X}_{i-\frac12,j}}{\Delta x}
+\frac{\mathcal{P}_{i,j+\frac12}-\mathcal{P}_{i,j-\frac12}}{\Delta p}.
\end{equation}
Here, $\mathcal{T}$ denotes the spatial transport flux associated with
$-(p_q/m)\partial_{x_q}W$,
$\mathcal{X}$ collects the spatial and cross-diffusion contributions
$2D_{xp}\partial_{x_qp_q}W+D_{xx}\partial_{x_qx_q}W$,
and $\mathcal{P}$ contains the momentum drift and diffusion contributions
$2\gamma\partial_{p_q}(p_qW)+D_{pp}\partial_{p_qp_q}W$.

For the transport term, the direction of propagation is determined by
$v_j=p_{q,j}/m$.
We use a second-order upwind reconstruction in the interior,
\begin{equation}
\mathcal{T}_{i+\frac12,j}
=
\begin{cases}
v_j(3W_{i,j}-W_{i-1,j})/2, & v_j\geq 0,\\
v_j(3W_{i+1,j}-W_{i+2,j})/2, & v_j<0,
\end{cases}
\end{equation}
together with the zero flux boundary condition.

The spatial diffusion and cross-diffusion contributions are combined into the flux
\begin{equation}
\mathcal{X}_{i+\frac12,j}
=
D_{xx}\frac{W_{i+1,j}-W_{i,j}}{\Delta x}
+
D_{xp}\left[(\delta_pW)_{i+1,j}+(\delta_pW)_{i,j}\right].
\end{equation}
The operator $\delta_p$ denotes the standard second-order centered difference in the interior,
\begin{equation}
(\delta_pW)_{i,j}
=
\frac{W_{i,j+1}-W_{i,j-1}}{2\Delta p},
\end{equation}
with second-order one-sided differences at the two momentum boundaries.
With this construction, the divergence of the cross-diffusion flux gives a second-order approximation of
$2D_{xp}\partial_{x_qp_q}W$.

The momentum drift and diffusion terms are treated together at the momentum interfaces
$p_{q,j+\frac12}=p_{\min}+j\Delta p$.
Using a centered approximation of the drift flux gives
\begin{equation}
\mathcal{P}_{i,j+\frac12}
=
\gamma p_{q,j+\frac12}(W_{i,j+1}+W_{i,j})
+
D_{pp}\frac{W_{i,j+1}-W_{i,j}}{\Delta p}.
\end{equation}
Thus, the drift contribution uses the arithmetic average of $W$ at the interface to approximate the flux
$2\gamma p_qW$.
Zero-flux boundary conditions are imposed for the diffusive and momentum fluxes at the truncated boundaries.

An important consequence of this construction is that each component of $L_1$ acts only on one or two physical modes.
Therefore, the discrete action of $L_1$ can be evaluated directly in TT form by modifying the corresponding cores and taking flux differences, without constructing the full discretization matrix.
In each substep of the projector splitting scheme, the contribution of $L_1[W]$ is projected to the subspace generated by the corresponding basis.
In our implementation, $L_1[W]$ is evaluated directly in TT form and then the projections becomes TT contractions between a tensor train and current TT basis.

We next introduce the Fourier grid used in the evaluation of $L_2$. 
The momentum discretization naturally determines the grid for the conjugate variable $\eta_q$, so no independent truncation of the $\eta_q$ domain is required.
The Fourier spacing is
\begin{equation}
\Delta\eta
=
\frac{2\pi}{N_p\Delta p}
=
\frac{2\pi}{p_{\max}-p_{\min}},
\end{equation}
and, for even $N_p$, the discrete Fourier frequencies are
\begin{equation}
\eta_{q,\ell}
=
\ell\Delta\eta,\qquad
\ell=0,1,\ldots,\frac{N_p}{2}-1,
-\frac{N_p}{2},\ldots,-1.
\end{equation}
Equivalently, after a Fourier shift, the grid covers
$[-\pi/\Delta p,\pi/\Delta p)$.

Because the momentum grid is cell-centered, the discrete Fourier transform contains a phase correction associated with the offset of the first momentum node.
Consistent with the Fourier convention used in the Wigner pseudo-differential operator, we approximate
\begin{equation}
\widehat{W}(\eta_{q,\ell})
\approx
\Delta p
\sum_{j=1}^{N_p}
W(p_{q,j})
\e^{-\ii p_{q,j}\eta_{q,\ell}},
\end{equation}
with the corresponding inverse transform
\begin{equation}
W(p_{q,j})
\approx
\frac{\Delta\eta}{2\pi}
\sum_\ell
\widehat{W}(\eta_{q,\ell})
\e^{\ii p_{q,j}\eta_{q,\ell}}.
\end{equation}
These transforms are evaluated using the fast Fourier transform (FFT) and inverse FFT together with the phase correction induced by the cell-centered momentum grid.

The same Fourier grid is used for the $\eta_q$ dependence of $\delta V(\bx,\be)$.
For a CP representation, the momentum-dependent factors
$D_\mu^{(2q)}(\eta_q)$
are evaluated on these Fourier frequencies.
For a TT representation, the corresponding even potential cores
$B_{\nu_{2q-1}\nu_{2q}}^{(2q)}(\eta_q)$
are discretized on the same grid.
With these definitions, the discrete representation of $L_2$ retains the same structural form as its continuous counterpart.

\section{Computational cost}
\label{sec_computational_cost}

In this section, we estimate the computational cost of one time step of the proposed method.
We separate our discussion to the cases $d\geq 2$ (typically 2 or 3) and $d=1$.
When the method is applied to $d=1$, the scaling of several operations is slightly different.

We first recall the notations used in the complexity analysis:
\begin{itemize}
    \item $d$: the physical dimension of the problem;
    \item $N_x$: the number of grid points in each spatial direction;
    \item $N_p$: the number of grid points in each momentum direction;
    \item $(r_1,r_2,\ldots,r_{2d-1})$: the prescribed TT ranks of the solution. For simplicity, we use $r$ as a uniform rank;
    \item $\varrho$: the number of terms in the CP representation of $\delta V$;
    \item $(\lambda_1,\ldots,\lambda_{2d-1})$: the TT ranks of $\delta V$ when the TT representation is used. For simplicity, we define $\lambda=\max_j\lambda_j$.
\end{itemize}
\subsection{Computational cost for $d\geq2$}
We first consider the case for $d\geq 2$ with a CP representation of $\delta V$.
The computational cost is estimated separately for the center-core and bond substeps.

We first consider a spatial center-core step for the CP representation.
The contribution of $F_2^{(k)}[W]$ is not explicitly assembled in our implementation.
Instead, the left and right environment matrices
$\mathsf L_\mu^{(k-1)}$
and
$\mathsf R_\mu^{(k)}$
are first constructed for each $\mu=1,\ldots,\varrho$, and these matrices are subsequently used in the matrix-free exponential action described in \cref{sec_matrix_free_implementation}.

According to
\cref{eq_L2_spatial_center_core_left_envorinment,eq_L2_spatial_center_core_right_envorinment},
the construction of the environments involves Fourier transforms in the momentum directions.
Using the fast Fourier transform, these operations require
$\mathcal{O}(d\varrho r^2N_p\log N_p)$
operations.
The elementwise multiplication of the TT cores by the corresponding CP factors requires
$\mathcal{O}(d\varrho r^2(N_x+N_p))$,
while the contractions with the left- and right-orthogonal TT bases require
$\mathcal{O}(d\varrho r^3(N_x+N_p))$.
We refer to \cite{wang2026dynamical} for a related analysis of these tensor contractions.
Therefore, the preparation of the environment matrices has cost
$\mathcal{O}(d\varrho r^3(N_x+N_p)+d\varrho r^2N_p\log N_p)$,
where lower-order terms have been omitted.

Once the environment matrices have been constructed, no additional summation over the CP terms is performed to explicitly form
$F_2^{(k)}[W]$
or
$M_k(x_q)$.
Instead, the summation over $\mu$ is carried out inside each matrix-vector product required by the exponential-action algorithm.
As discussed in \cref{sec_matrix_free_implementation}, one application of $M_k(x_q)$ costs
$\mathcal{O}(\varrho r^3)$.
Since the exponential action is evaluated independently at each of the $N_x$ spatial grid points, and requires $n_{\mathrm{mv}}$ matrix-vector products, the corresponding cost is
$\mathcal{O}(n_{\mathrm{mv}}\varrho r^3N_x)$.
Here and below, $n_{\mathrm{mv}}$ denotes a representative number of matrix-vector products required for one local exponential action.

The construction of $F_1^{(k)}[W]$ does not involve the CP rank $\varrho$.
As discussed in \cref{sec_discretization}, the differential operators are applied directly to the TT representation, followed by projection onto the current TT basis.
The corresponding contractions have cost of order
$\mathcal{O}(dr^3(N_x+N_p))$,
which is typically lower than the cost associated with the potential term when $\varrho$ is larger than one.

Combining these contributions, one spatial center-core step requires
\begin{equation}
\mathcal{O}\left(
d\varrho r^3(N_x+N_p)
+
d\varrho r^2N_p\log N_p
+
n_{\mathrm{mv}}\varrho r^3N_x
\right)
\end{equation}
operations.

A momentum center-core step can be analyzed in the same way.
The environment matrices are constructed with the same leading-order cost, while the exponential action is now evaluated at each of the $N_p$ momentum grid points.
Therefore, one momentum center-core step requires
\begin{equation}
\mathcal{O}\left(
d\varrho r^3(N_x+N_p)
+
d\varrho r^2N_p\log N_p
+
n_{\mathrm{mv}}\varrho r^3N_p
\right)
\end{equation}
operations.

The bond substeps are less expensive because the exponential action is applied only to an $r\times r$ bond matrix and is not repeated over a physical grid.
The corresponding matrix-free exponential therefore requires
$\mathcal{O}(n_{\mathrm{mv}}\varrho r^3)$
operations, in addition to the environment contractions required for the bond update.

A complete left-to-right projector-splitting sweep contains $2d$ center-core steps and $2d-1$ bond steps.
Including the orthogonalization operations required to move the TT center during the sweep, the leading computational cost of one time step for a CP representation of $\delta V$ can therefore be estimated by
\begin{equation}
\label{eq_computational_cost}
\mathcal{O}\!\left(
d^2\varrho r^3(N_x+N_p)
+
d^2\varrho r^2N_p\log N_p
+
d\,n_{\mathrm{mv}}\varrho r^3(N_x+N_p)
\right).
\end{equation}
For fixed physical dimension $d=2$ or $3$, this may be summarized as linear scaling with $N_x$ and $N_p$, up to the FFT factor $\log N_p$, and cubic scaling with the prescribed solution TT rank $r$ when $n_{\mathrm{mv}}$ and $\varrho$ remain moderate.

We next consider the case in which the potential difference is represented in TT form.
The same operations are required, but the global CP summation over $\mu$ is replaced by contractions over the neighboring potential TT indices.
In particular, the update of an environment across a potential core involves at most
$\lambda_{j-1}\lambda_j$
combinations of potential bond indices.
Under the uniform bound $\lambda_j\leq\lambda$, the preparation of the potential environments for a center-core step therefore costs
\begin{equation}
\mathcal{O}\left(
d\lambda^2r^3(N_x+N_p)
+
d\lambda^2r^2N_p\log N_p
\right).
\end{equation}

For a spatial or momentum center-core step, one application of the local matrix involves the two neighboring potential TT indices $\nu_{k-1}$ and $\nu_k$.
As discussed in \cref{sec_matrix_free_implementation}, one matrix-vector product therefore costs
$\mathcal{O}(\lambda_{k-1}\lambda_k r^3)$,
which is bounded by
$\mathcal{O}(\lambda^2r^3)$.
The matrix-free exponential actions over the spatial and momentum grids consequently cost
$\mathcal{O}(n_{\mathrm{mv}}\lambda^2r^3N_x)$
and
$\mathcal{O}(n_{\mathrm{mv}}\lambda^2r^3N_p)$,
respectively.

For a bond step, only the potential bond index $\nu_k$ remains in
\cref{eq_bond_G2_TT}.
Thus, one matrix-vector product costs
$\mathcal{O}(\lambda_k r^3)$,
and the corresponding exponential action costs
$\mathcal{O}(n_{\mathrm{mv}}\lambda r^3)$
under the uniform-rank estimate.

Combining the center-core steps, bond steps, environment construction, and orthogonalization, the leading cost of one time step with a TT representation of $\delta V$ can be bounded by
\begin{equation}
\label{eq_computational_cost_TT_potential}
\mathcal{O}\!\left(
d^2\lambda^2r^3(N_x+N_p)
+
d^2\lambda^2r^2N_p\log N_p
+
d\,n_{\mathrm{mv}}\lambda^2r^3(N_x+N_p)
\right),
\end{equation}
where the lower-order bond exponential contribution
$\mathcal{O}(d\,n_{\mathrm{mv}}\lambda r^3)$
has been absorbed into the leading-order estimate.
Therefore, the CP and TT representations of the potential difference lead to the same dependence on the solution rank and grid resolution.
The difference lies in the complexity of the potential representation: the CP formulation depends linearly on the number $\varrho$ of separable terms, whereas the TT formulation depends locally on the neighboring potential TT ranks and is bounded by a quadratic dependence on the uniform potential rank $\lambda$.

\subsection{Computational cost for $d=1$}
For the special case $d=1$, the tensor train for the solution contains only two physical modes, $x_1$ and $p_1$, and the dimensions of the center-core subproblems are smaller than those for $d=2,3$.
In particular, the spatial center core $C^{(1)}(x)$ has size $1\times r$, while the momentum center core $C^{(2)}(p)$ has size $r\times 1$.
Consequently, the local matrices appearing in the two center-core exponential updates have dimension $r\times r$, rather than $r^2\times r^2$.

For a CP representation of $\delta V$ with $\varrho$ terms, one matrix-vector product in the matrix-free exponential action therefore requires
$\mathcal{O}(\varrho r^2)$
operations.
The exponential-action costs of the spatial and momentum center-core steps are thus
$\mathcal{O}(n_{\mathrm{mv}}\varrho r^2N_x)$
and
$\mathcal{O}(n_{\mathrm{mv}}\varrho r^2N_p)$,
respectively.
The bond matrix $S^{(1)}$, however, has size $r\times r$.
The corresponding vectorized bond equation therefore still involves an operator acting on a vector of dimension $r^2$, and one matrix-free operator application costs
$\mathcal{O}(\varrho r^3)$.
Since the bond exponential action is performed only once per sweep and is not repeated over the physical grids, its total cost is
$\mathcal{O}(n_{\mathrm{mv}}\varrho r^3)$.

Combining the center-core updates, the bond update, Fourier transforms, and the required contractions, the dominant cost of one time step for $d=1$ with a CP representation of the potential difference can be estimated as
\begin{equation}
\mathcal{O}\!\left(
n_{\mathrm{mv}}\varrho r^2(N_x+N_p)
+
n_{\mathrm{mv}}\varrho r^3
+
\varrho r^2(N_x+N_p)
+
\varrho rN_p\log N_p
\right).
\end{equation}
For $d=1$, a TT representation of $\delta V(x,\eta)$ contains only one potential bond rank $\lambda_1$ and is equivalent in structure to a separable representation with $\lambda_1$ terms.
Therefore, the same one-dimensional complexity estimate applies with $\varrho$ replaced by $\lambda_1$.

Hence, the one-dimensional case is substantially cheaper than the $d=2,3$ cases, not only because fewer TT modes are involved, but also because the center-core exponential operators have dimension $r\times r$ rather than $r^2\times r^2$.

\section{Numerical examples}
\label{sec_Numerical_Examples}
In this section, we present several numerical examples to validate our method.
Since the contribution from $L_2$ is advanced through the exponential operator, it does not introduce an additional explicit stability restriction. 
Therefore, in all numerical experiments, the time step is chosen according to a CFL-type stability condition associated with the explicitly treated operator $L_1$. 
In the implementation, we estimate
\begin{equation}
\label{dt_timestep}
\Delta t
=\mathrm{CFL}
\left(\sum_{q=1}^{d}
\left(
\frac{3}{2}\frac{\max_j|p_{q,j}|}{m\Delta x}
+\frac{2D_{xx}}{\Delta x^2}
+\frac{2D_{pp}}{\Delta p^2}
+\frac{2|D_{xp}|}{\Delta x\Delta p}
+2|\gamma|
\right)\right)^{-1}
\end{equation}
where the factor $3/2$ accounts for the second-order upwind discretization of the spatial transport term. The time step used in the simulations is chosen not to exceed this estimate. This condition is intended as a practical stability bound for the forward-Euler treatment of $L_1$, rather than a sharp stability condition for the full projector-splitting scheme.
All the following numerical experiments were conducted on a MacBook Pro equipped with an Apple M4 CPU.

The numerical examples in this section are designed to examine the proposed method from both numerical and physical perspectives.
We first validate the method against existing one- and two-dimensional results.
We then study how environmental friction and diffusion affect the phase-space dynamics and its low-rank structure. 
Finally, we demonstrate the feasibility of simulating a six-dimensional WFP problem.

\subsection{1D1P Wigner--Fokker--Planck equation}
For $d=1$, the Wigner function $W$ is represented by \cref{eq_Wigner_1D1V_three_forms}, and the proposed method reduces to a projector-splitting method for a bivariate function.
We follow the parameter setting of Example 4 in \cite{yi2026time}.
The coefficients in \cref{eq_wfp} are given by
\begin{equation}
m=1,\qquad
\hbar=0.1,\quad
\gamma=1,\qquad
D_{pp}=0.1,\qquad D_{xp}=0,\qquad D_{xx}=0.1,\qquad \mathrm{CFL}=0.1.
\end{equation}
These coefficients satisfy the Lindblad condition \cref{eq_lindblad_condition}. 
Therefore, the corresponding continuous density-operator evolution is completely positive.
Following the setting in \cite{yi2026time}, the initial Wigner function is chosen as
\begin{equation}
W(0,x,p)
=
\frac{1}{\pi\hbar}
\exp\left(
-\frac{(x-0.1)^2+(p+0.2)^2}{\hbar}
\right),
\end{equation}
and the potential is given by
\begin{equation}
V(x)=\arctan(10x)+\frac{\pi}{2}.
\end{equation}

Following the setting of the reference paper, both the spatial and momentum domains are truncated to $[-4,4]$, and we use $N_x=N_p=128$ grid points in the simulation.
Based on \cref{dt_timestep}, the time step is $\dt\approx 5.009 \times 10^{-4}$.
The rank of the numerical solution is fixed at $r=25$.

To construct the CP representation \cref{eq_deltaV_CP} of the potential difference $\delta V$, we first evaluate
\begin{equation}
\delta V(x_i,\eta_\ell)
=
V\left(x_i+\frac{\hbar}{2}\eta_\ell\right)
-
V\left(x_i-\frac{\hbar}{2}\eta_\ell\right)
\end{equation}
on the discrete $(x,\eta)$ grid.
Since $d=1$, $\delta V$ is represented by a matrix, and its singular value decomposition directly provides a separable representation in $x$ and $\eta$.
Singular values smaller than $10^{-10}$ are discarded, resulting in a rank-$47$ approximation of $\delta V$.

\begin{figure}
    \centering
    \includegraphics[width=0.8\linewidth]{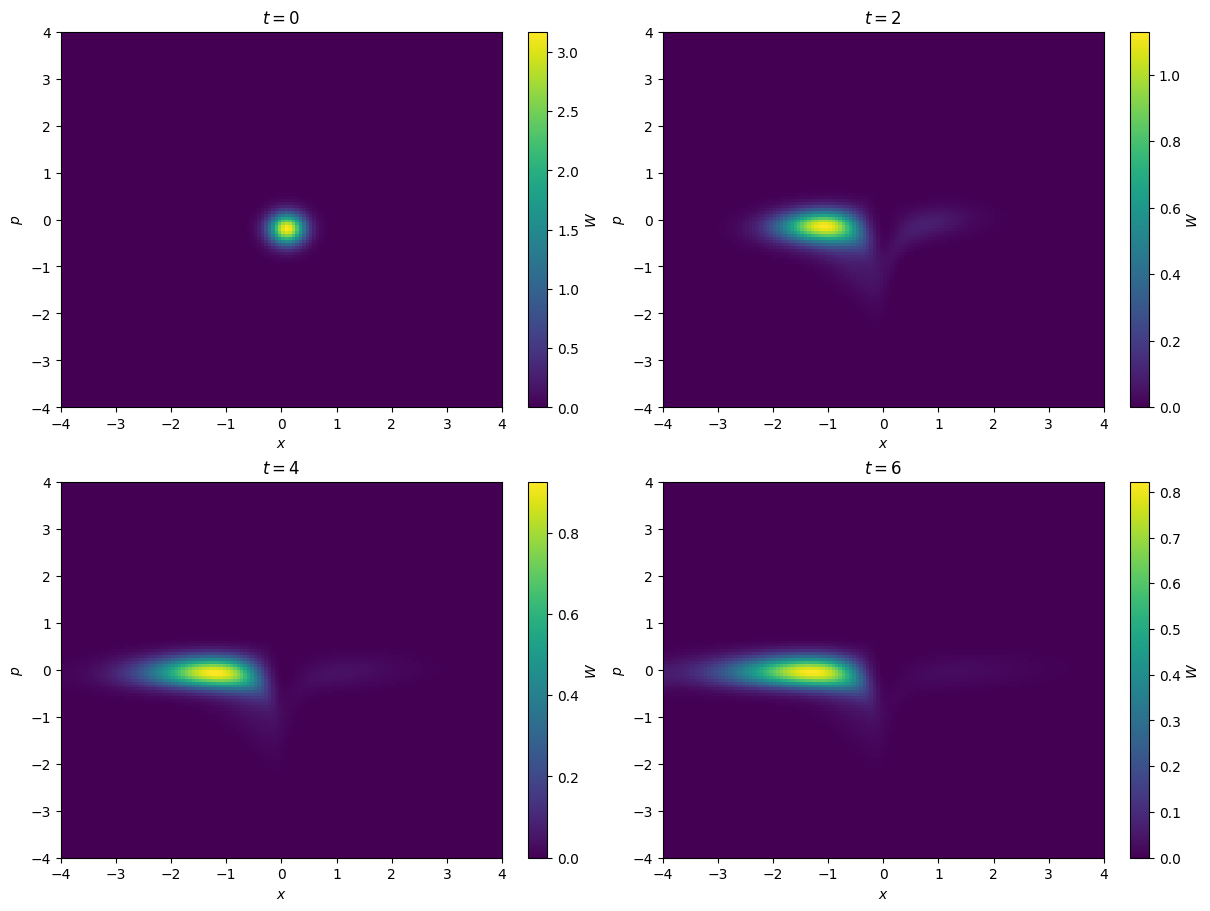}
    \caption{Wigner function evolution for potential $V(x) = \arctan(10x) + \frac{\pi}{2}$.}
    \label{fig_example_1D1V}
\end{figure}

We show the phase-space plots of the numerical solution at $t=0,2,4,$ and $6$ in \cref{fig_example_1D1V}.
The numerical results are consistent with those reported in \cite{yi2026time}, providing a validation of the proposed low-rank projector-splitting method.
The test also validates the method in a regime where coherent Wigner transport and environmental friction/diffusion are simultaneously present.

\subsection{2D2P Wigner equation}
The proposed method is designed in particular for problems with $d=2$ and $d=3$. In this section, we test the method in a two-dimensional physical space, corresponding to a four-dimensional position--momentum phase space.

Since numerical results for the Wigner--Fokker--Planck equation in two spatial dimensions are still limited in the existing literature, we follow the parameter setting of Example 4.3 in \cite{shao2026dynamical}. 
In this example, the Fokker--Planck terms are absent, and the coefficients in \cref{eq_wfp} are given by
\begin{equation}
m=1,\qquad \hbar=1,\qquad
\gamma=D_{pp}=D_{xp}=D_{xx}=0,\qquad \mathrm{CFL}=0.1.
\end{equation}
Although this example does not describe an open quantum system, it provides a benchmark against an existing dynamical low-rank method of Wigner function evolution.

The initial condition is
\begin{equation}
\label{eq_initial_fermion}
\begin{split}
W\left(0,x_1,x_2,p_1,p_2\right)
={}& \frac{1}{2\pi^2}\exp\left[-\frac{(x_1-x_1^0)^2}{2a^2}-\frac{(x_2-x_2^0)^2}{2a^2}-2a^2(p_1-p_1^0)^2-2a^2(p_2-p_2^0)^2\right] \\
&+\frac{1}{2\pi^2}\exp\left[-\frac{(x_1-x_2^0)^2}{2a^2}-\frac{(x_2-x_1^0)^2}{2a^2}-2a^2(p_1-p_2^0)^2-2a^2(p_2-p_1^0)^2\right] \\
&-\frac{1}{\pi^2}\exp\left[-\frac{(x_1-x_1^0)^2+(x_1-x_2^0)^2+(x_2-x_1^0)^2+(x_2-x_2^0)^2}{4a^2}\right] \\
&\times \exp\left[\frac{(x_1^0-x_2^0)^2}{4a^2}-2a^2\left(p_1-\frac{p_1^0+p_2^0}{2}\right)^2-2a^2\left(p_2-\frac{p_1^0+p_2^0}{2}\right)^2\right] \\
&\times \cos\left[(x_1^0-x_2^0)(p_1-p_2)-(p_1^0-p_2^0)(x_1-x_2)\right],
\end{split}
\end{equation}
and the interaction potential is
\begin{equation}
\label{eq_electron_potential}
V(x_1,x_2)=\frac{1}{\sqrt{(x_1-x_2)^2+1}}.
\end{equation}
In \cref{eq_initial_fermion}, the parameters are
$p_1^0=0.5$, $p_2^0=-0.5$, $x_1^0=-2$, $x_2^0=2$, and $a=2$.
The computational domains are truncated to
$[x_{\min},x_{\max}]^2=[-15,15]^2$ and
$[p_{\min},p_{\max}]^2=[-5\pi/6,5\pi/6]^2$,
with $N_x=150$ and $N_p=100$.
According to \cref{dt_timestep}, the time step in this numerical example is $\dt \approx 2.570\times 10^{-3}$.

The purpose of this example is to compare the numerical results obtained with our method with those reported in the existing literature.
However, since the Fokker--Planck terms are absent in this example, the dissipative mechanism that may enhance the low-rank structure of the solution is also absent.
Consequently, the solution may require relatively large TT ranks to achieve an accurate approximation.

For the potential given by \cref{eq_electron_potential}, the potential difference $\delta V$ does not immediately admit an explicit CP or TT decomposition.
In practice, when a direct separated representation is unavailable, sampling-based methods can be used to construct low-rank approximations without forming the full tensor.
For the TT format, TT-cross \cite{oseledets2010ttcross,savostyanov2014quasioptimality} is a well-developed approach that constructs a tensor-train approximation using adaptively selected tensor entries without explicitly forming the full tensor.
For the CP format, several sampling-based approaches have also been developed.
We refer the reader to \cite{bebendorf2013equi,battaglino2018practical,sorensen2019fiber} for representative methods for constructing CP approximations from sampled entries or fibers.

In the present example, however, we do not use sampling-based methods to construct the CP or TT representations of $\delta V$.
Instead, we construct an approximate CP decomposition through a truncated Fourier expansion of the one-dimensional kernel
\begin{equation}
v(z)=\frac{1}{\sqrt{z^2+1}}.
\end{equation}
We truncate $z$ to $[-L,L]$, extend $v$ periodically, and approximate it by the Fourier series
\begin{equation}
v(z)\approx c_0+\sum_{m=1}^{M}
\left(
c_m\e^{\ii\omega_m z}
+\overline{c_m}\e^{-\ii\omega_m z}
\right),
\qquad
\omega_m=\frac{\pi m}{L}.
\end{equation}
In the numerical implementation, the Fourier coefficients are approximated by the discrete Fourier transform
\begin{equation}
c_m\approx
\frac{1}{N_s}\sum_{j=0}^{N_s-1}
v(z_j)\e^{-2\pi\ii m j/N_s}\e^{\ii\omega_m L},
\qquad
z_j=-L+\frac{2Lj}{N_s},
\end{equation}
where $j=0,\ldots,N_s-1$.
We choose $L=100$ and $N_s=1000$.

Since the constant Fourier mode cancels in the potential difference, the remaining contribution can be written as
\begin{equation}
\delta V
\approx
\sum_{m=1}^{M}\sum_{\sigma=\pm1}T_{m,\sigma}
+
\overline{
\sum_{m=1}^{M}\sum_{\sigma=\pm1}T_{m,\sigma}
},
\end{equation}
where
\begin{equation}
T_{m,\sigma}
=
\left(\sigma c_m\e^{\ii\omega_m x_1}\right)
\left(\e^{\ii\sigma\hbar\omega_m\eta_1/2}\right)
\left(\e^{-\ii\omega_m x_2}\right)
\left(\e^{-\ii\sigma\hbar\omega_m\eta_2/2}\right).
\end{equation}
Each $T_{m,\sigma}$ is a rank-one separable function of
$(x_1,\eta_1,x_2,\eta_2)$.
Therefore, when the complex-conjugate terms are expanded explicitly, the resulting CP representation contains $4M$ rank-one terms.

For comparison, we also construct a TT approximation of $\delta V$.
In this test, rather than applying TT-cross, we first construct the full discrete tensor of size
$N_x\times N_p\times N_x\times N_p$
and then apply the standard TT-SVD algorithm \cite{oseledets2011tensor}.

We consider different TT ranks and different numbers of CP terms in our setting.
To evaluate the approximation accuracy, we reconstruct the corresponding full tensors and compute their relative $L^2$ errors with respect to the directly evaluated tensor.
The results are reported in \cref{table_deltaV_error}.
\begin{table}[t]
    \centering
    \begin{tabular}{c|c|c|c|c}
        \hline
        TT rank $(\lambda_1,\lambda_2,\lambda_3)$ &
        Relative $L^2$ error &
        CP terms $\varrho$ &
        Relative $L^2$ error \\
        \hline
        $(40,196,71)$ &
        $9.54\times10^{-3}$ &
        $448$ &
        $9.95\times10^{-3}$ \\
        $(102,558,100)$ &
        $7.98\times10^{-6}$ &
        $1264$ &
        $9.91\times10^{-6}$ \\
        $(150,879,100)$ &
        $5.10\times10^{-9}$ &
        $1996$ &
        $2.21\times10^{-15}$ \\
        $(150,888,100)$ &
        $1.25\times10^{-13}$ &
        $1996$ &
        $2.21\times10^{-15}$ \\
        \hline
    \end{tabular}
    \caption{TT ranks and numbers of CP terms required for different approximation tolerances of the potential difference $\delta V$, together with the corresponding relative $L^2$ errors.}
    \label{table_deltaV_error}
\end{table}

The TT approximations in \cref{table_deltaV_error} are constructed by TT-SVD, which provides a quasi-optimal approximation in the Frobenius norm for prescribed TT ranks \cite{oseledets2011tensor}.
In contrast, the CP approximations used here are obtained from the particular Fourier construction above and are not expected to be optimal among all CP representations of the same size.
Nevertheless, even with this specialized construction, the number of CP terms $\varrho$ is substantially smaller than the products of neighboring TT ranks $\lambda_{k-1}\lambda_k$ appearing in the corresponding center-core operations.
According to the computational cost analysis in the previous section, the CP representation is therefore more efficient for this particular potential and implementation.
We consequently use the CP representation with $\varrho=448$ terms in the remainder of this example.

\begin{figure}
    \centering
    \includegraphics[width=\linewidth]{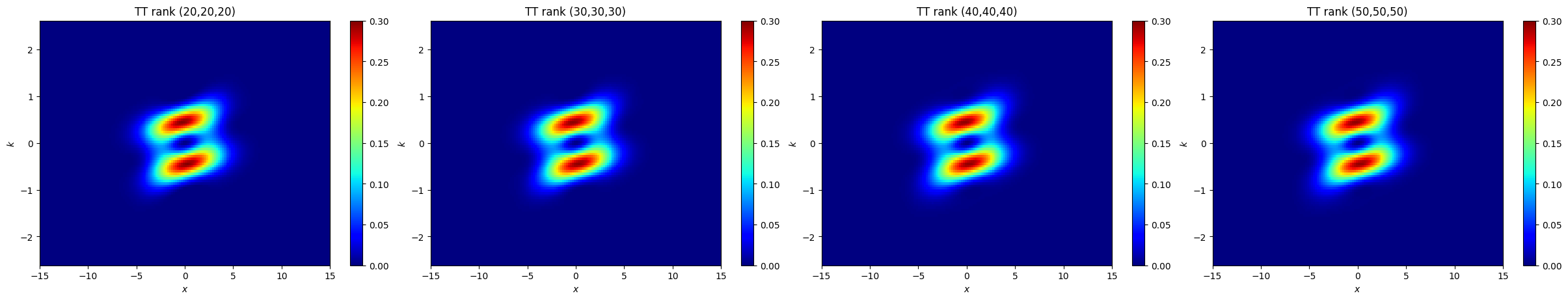}
    \caption{$W_1(x,p)+W_2(x,p)$ for different TT ranks at $t=3$.}
    \label{fig_convergence}
\end{figure}

We next investigate the dependence of the numerical solution on the prescribed TT rank of $W$.
We consider the TT ranks
$(10,10,10)$, $(20,20,20)$, $(30,30,30)$, $(40,40,40)$, and $(50,50,50)$,
and evolve the system up to $t=3$.
The corresponding phase-space quantity $W_1(x,p)+W_2(x,p)$ is shown in \cref{fig_convergence}.
Here, $W_1$ and $W_2$ denote the reduced one-particle Wigner functions, defined by
\begin{equation}
W_1(t,x,p)
=
\iint_{\mathbb{R}\times\mathbb{R}}
W(t,x,x_2,p,p_2)\,\dd x_2\,\dd p_2,
\qquad
W_2(t,x,p)
=
\iint_{\mathbb{R}\times\mathbb{R}}
W(t,x_1,x,p_1,p)\,\dd x_1\,\dd p_1.
\end{equation}
In addition, using the result obtained with TT rank $(50,50,50)$ as a reference, we report the maximum-norm and $L^2$-norm differences for the lower-rank solutions in \cref{table_convergence}.

\begin{table}[t]
    \centering
    \begin{tabular}{c|c|c}
        TT rank & maximum difference & $L^2$ difference \\ \hline
        $(10,10,10)$ & 0.009905 & 0.020681 \\ \hline
        $(20,20,20)$ & 0.007332 & 0.015996 \\ \hline
        $(30,30,30)$ & 0.004224 & 0.009041\\ \hline
        $(40,40,40)$ & 0.000468 & 0.000889
    \end{tabular}
    \caption{Differences between the results obtained with different TT ranks and the reference result with TT rank $(50,50,50)$.}
    \label{table_convergence}
\end{table}

\begin{figure}
    \centering
    \includegraphics[width=\linewidth]{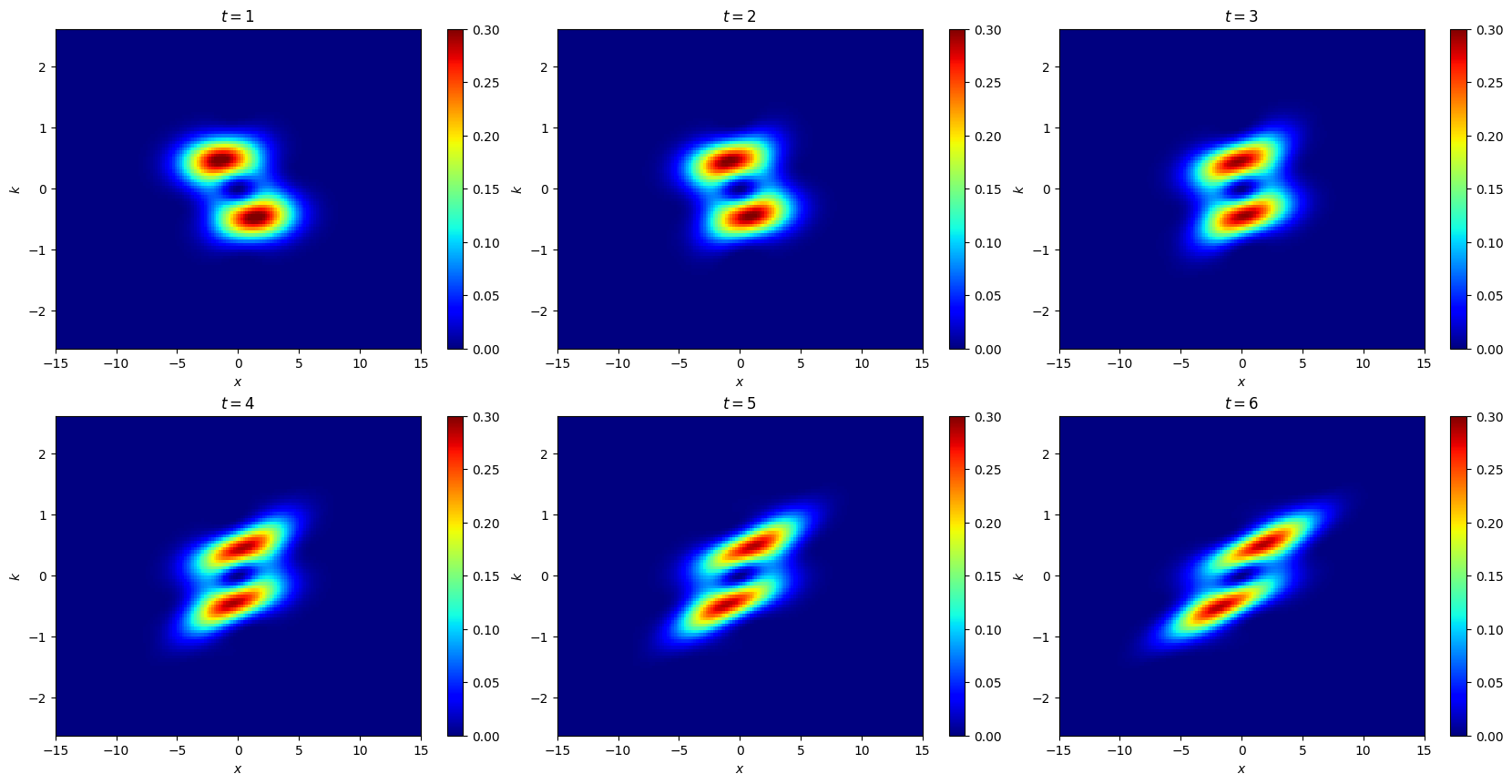}
    \caption{Evolution of $W_1(x,p)+W_2(x,p)$ with TT rank $(30,30,30)$.}
    \label{fig_example_Wigner_2D2V}
\end{figure}

As shown in \cref{fig_convergence}, the phase-space plots obtained with different TT ranks are visually very similar.
The quantitative differences reported in \cref{table_convergence} decrease consistently as the prescribed TT rank increases, indicating convergence with respect to the TT rank.
Since the solutions obtained with moderate ranks already reproduce the main qualitative features of the reference solution, we choose the TT rank $(30,30,30)$ for the longer-time simulation up to $t=6$.
The resulting evolution is shown in \cref{fig_example_Wigner_2D2V}.
The evolution is qualitatively consistent with the corresponding results reported in \cite{shao2026dynamical}.

Having validated the low-rank treatment of the coherent Wigner dynamics in this four-dimensional phase-space setting, we next turn to the effect of environmental friction and diffusion on both the physical dynamics and the numerical low-rank structure.

\subsection{2D2P open-system dynamics with a weakly nonquadratic potential}
We next test the proposed method for a two-dimensional quadratic potential with a weak trigonometric perturbation,
\begin{equation}
\label{eq_perturbed_quadratic}
    V(x_1,x_2)
    =
    \frac{1}{2}(x_1^2+x_2^2)
    +x_1+x_2
    +0.1\sin(x_1+x_2).
\end{equation}
This potential can be viewed as a natural two-dimensional extension of an example considered in \cite{yi2026time}.
The initial condition is chosen as a weighted combination of two Gaussian wave packets,
\begin{equation}
\begin{split}
    W(0,\bx,\bp)
    &=
    \frac{3}{5(\pi\hbar)^2}
    \exp\left(
    -\frac{(x_1-0.5)^2+(x_2+0.4)^2+(p_1+0.6)^2+(p_2-0.3)^2}{\hbar}
    \right) \\
    &\quad+
    \frac{2}{5(\pi\hbar)^2}
    \exp\left(
    -\frac{(x_1+0.7)^2+(x_2-0.6)^2+(p_1-0.5)^2+(p_2+0.4)^2}{\hbar}
    \right).
\end{split}
\end{equation}
We set $\hbar=0.1$ and the particle mass $m=1$.

To examine the effect of the Fokker--Planck terms on the dynamics, we consider three sets of parameters:
\begin{enumerate}[(a)]
    \item $\gamma=0$, $D_{pp}=0$, $D_{xp}=0$, $D_{xx}=0$;
    \item $\gamma=1$, $D_{pp}=0.1$, $D_{xp}=0$, $D_{xx}=0.1$;
    \item $\gamma=2$, $D_{pp}=0.2$, $D_{xp}=0$, $D_{xx}=0.2$.
\end{enumerate}
Case (a) contains no Fokker--Planck contribution and therefore reduces to the Wigner equation.
It describes the dynamics of a closed quantum system.
Cases (b) and (c) introduce friction and diffusion of increasing strength, allowing us to compare the purely Wigner evolution with increasingly dissipative Wigner--Fokker--Planck dynamics.
Both parameter sets satisfy the Lindblad condition \cref{eq_lindblad_condition}, so the corresponding continuous density-operator evolutions are completely positive.

For the potential in \cref{eq_perturbed_quadratic}, the difference potential admits an exact CP representation with six separable terms. In particular,
\begin{equation}
\begin{split}
    \delta V(\bx,\be)
    &=
    \hbar(1+x_1)\eta_1
    +\hbar(1+x_2)\eta_2
    +0.2\cos(x_1)\cos(x_2)
    \sin\left(\frac{\hbar\eta_1}{2}\right)
    \cos\left(\frac{\hbar\eta_2}{2}\right) \\
    &\quad
    +0.2\cos(x_1)\cos(x_2)
    \cos\left(\frac{\hbar\eta_1}{2}\right)
    \sin\left(\frac{\hbar\eta_2}{2}\right) \\
    &\quad
    -0.2\sin(x_1)\sin(x_2)
    \cos\left(\frac{\hbar\eta_1}{2}\right)
    \sin\left(\frac{\hbar\eta_2}{2}\right) \\
    &\quad
    -0.2\sin(x_1)\sin(x_2)
    \sin\left(\frac{\hbar\eta_1}{2}\right)
    \cos\left(\frac{\hbar\eta_2}{2}\right).
\end{split}
\end{equation}
Thus, no additional low-rank approximation of $\delta V$ is required in this example.

We truncate both the spatial and momentum domains to $[-4,4]^2$.
For each spatial direction, we use $N_x=150$ grid points, while each momentum direction is discretized using $N_p=100$ grid points.
The prescribed TT rank is $(40,40,40)$.
Under this setting, the TT representation of $W$ requires $410000$ stored entries, compared with $225000000$ entries for the corresponding full tensor representation.

\begin{figure}
    \centering
    \includegraphics[width=\linewidth]{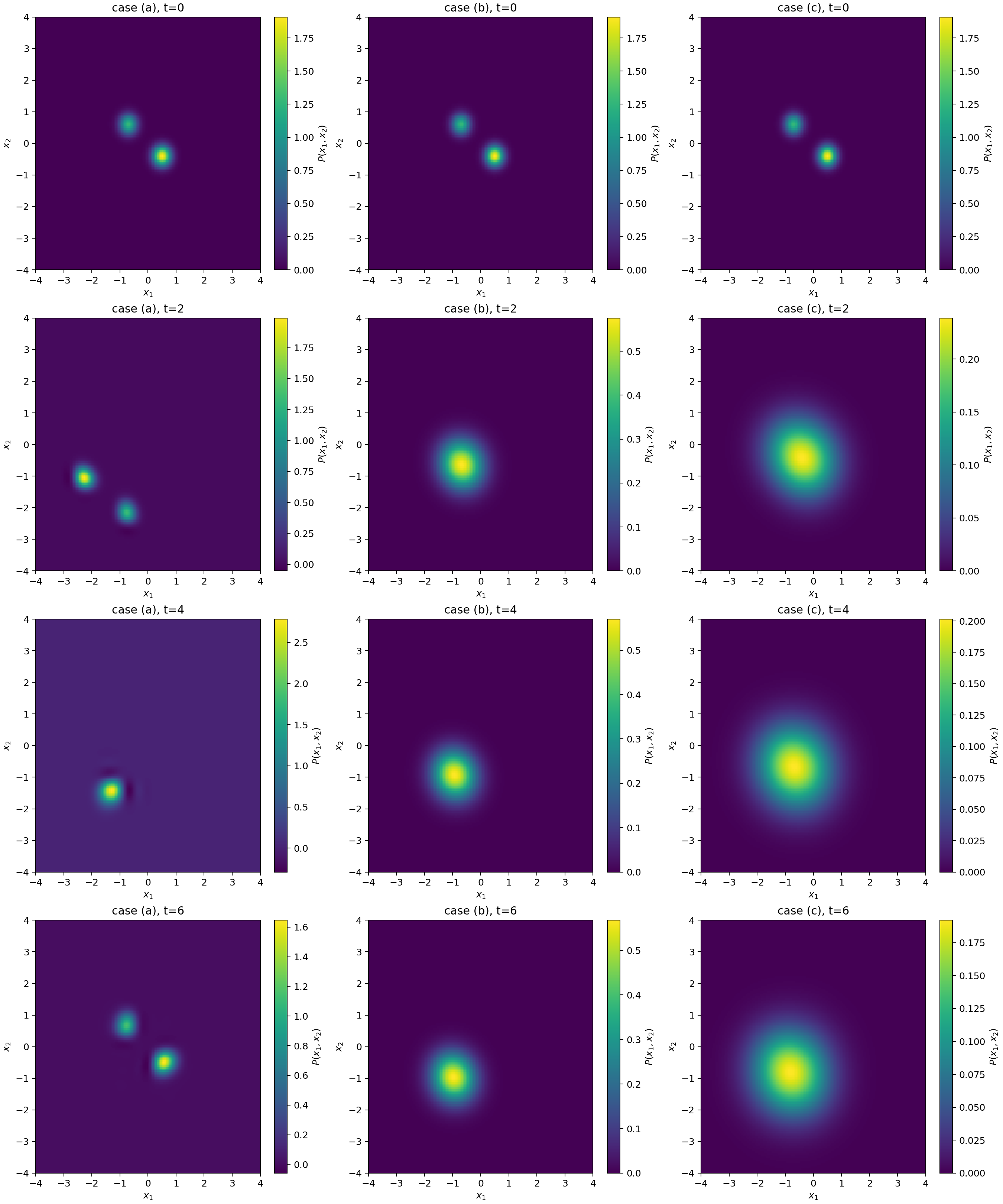}
    \caption{Position density at $t=0,2,4,6$ for the three Fokker--Planck parameter settings.}
    \label{fig_example_2D2V_different_FP}
\end{figure}

The system is evolved up to $t=6$ under the three parameter settings.
To visualize the dynamics, we consider the position density
\begin{equation}
    P(t,\bx)
    =
    \iint_{\mathbb{R}^2}
    W(t,\bx,\bp)\,\dd\bp,
\end{equation}
which is evaluated numerically over the truncated momentum domain.
The position densities at $t=0,2,4,6$ are shown in \cref{fig_example_2D2V_different_FP}.

As shown in \cref{fig_example_2D2V_different_FP}, the Fokker--Planck terms have a clear dissipative effect on the evolution.
In the absence of the Fokker--Planck contribution, the two wave packets retain more pronounced dynamical structures over time.
As the friction and diffusion coefficients increase, these structures are progressively damped, and the two initially separated wave packets gradually merge and concentrate around the low-potential region.
The damping is stronger in case (c) than in case (b), consistent with the larger friction and diffusion coefficients.
From the open-system perspective, this comparison illustrates the competition between coherent quantum evolution and environment-induced relaxation.
Case (a) contains only the Hamiltonian Wigner dynamics, whereas cases (b) and (c) additionally include friction and diffusion generated by the environment.
Over the simulated time interval, these environmental terms damp and smooth the spatial structures visible in the position density generated by the coherent dynamics, with a stronger effect observed for the larger friction and diffusion coefficients in case (c).

\begin{figure}
    \centering
    \includegraphics[width=\linewidth]{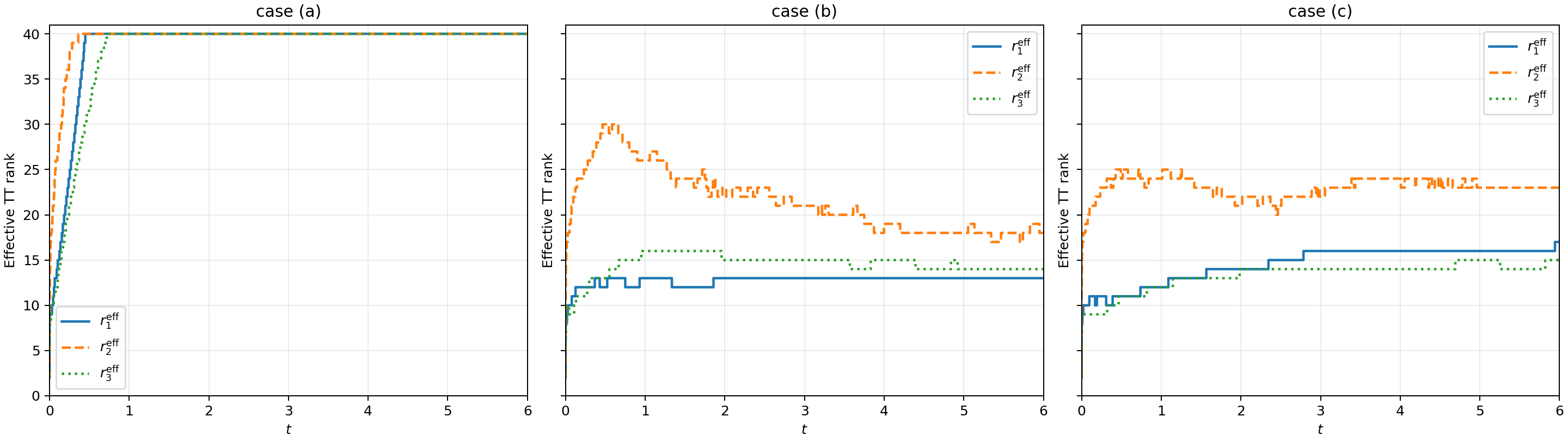}
    \caption{Evolution of the effective TT ranks with $\delta=10^{-6}$.}
    \label{fig_example_2D2V_effective_ranks}
\end{figure}

We next examine how this dissipative behavior affects the numerical low-rank structure of the solution.
Following \cite{wang2026dynamical,wang2026implicit}, we define the effective rank associated with the $k$th TT bond by
\begin{equation}
    r_k^{\mathrm{eff}}
    =
    \max\left\{
    \ell:\sigma_\ell\geq \delta\sigma_1
    \right\},
\end{equation}
where $\sigma_1\geq\cdots\geq\sigma_r$ are the singular values of the bond matrix $S^{(k)}$ when the solution is written in the bond form in \cref{eq_bond_form}.
For $d=2$, the TT representation contains three bond matrices, corresponding to $k=1,2,3$, and therefore three effective TT ranks.
We choose $\delta=10^{-6}$ and plot their evolution in \cref{fig_example_2D2V_effective_ranks}.

The evolution of the effective TT ranks shows that the introduction of the Fokker--Planck terms leads to a substantial reduction in the numerical ranks required to represent the solution.
A similar reduction of numerical complexity due to collisional or dissipative effects has also been observed for low-rank approximations of kinetic equations \cite{wang2026dynamical,wang2026implicit}.
For the present Wigner--Fokker--Planck example, the result suggests that dissipation and diffusion suppress some of the fine and oscillatory phase-space structures generated by the unitary Wigner dynamics, thereby increasing the compressibility of the Wigner function.
This behavior is also qualitatively consistent with the physical picture of decoherence and dissipation, in which interaction with an environment tends to suppress coherent quantum structures.
While the present experiment does not establish that open quantum dynamics is generally of lower rank than the corresponding closed-system dynamics, it indicates that environmental effects can, in some regimes, make the reduced phase-space dynamics more amenable to low-rank approximation.
From a numerical perspective, this observation shows that the environment can affect not only the physical evolution but also the compressibility of the solution. In the present example, the dissipative mechanisms that smooth the phase-space dynamics simultaneously reduce the effective tensor ranks required to represent it.

\subsection{3D3P simulation of an open quantum system}
To further demonstrate the capability of the proposed method in higher dimensions, we consider a three-dimensional example with $d=3$, corresponding to a six-dimensional position--momentum phase space.
This example is motivated by \cite{arnold2012wigner}, where the authors analyze Wigner–Fokker–Planck equations with smooth perturbations of a quadratic potential,
\begin{equation}
V(\bx)
=
\frac{1}{2}|\bx|^2+\lambda V_0(\bx),
\end{equation}
The choice $
V_0(\bx)=\sin(\bx\cdot\mathbf{k}_0)$ is also discussed there as an illustrative example.
In the following numerical experiment, we choose $\lambda=0.2$ and $\mathbf{k}_0=(1,1,1)^\top$.
In this setting, the potential difference can be written as
\begin{equation}
    \delta V(\bx,\be)
    =
    \hbar\sum_{j=1}^3 x_j\eta_j
    +
    \frac{\lambda}{2\ii}
    \sum_{\sigma=\pm1}
    \sum_{\tau=\pm1}
    \tau
    \prod_{j=1}^3
    \e^{\ii\sigma x_j}
    \e^{\ii\tau\hbar\eta_j/2}.
\end{equation}
This expression provides an exact CP representation of $\delta V$ with seven separable terms: three terms from the quadratic potential and four terms from the trigonometric perturbation.

The initial Wigner function is taken to be a displaced Gaussian,
\begin{equation}
W(0,\bx,\bp)
=
\frac{1}{(\pi\hbar)^3}
\exp\left(
-\frac{|\bx-\bx_0|^2+|\bp-\bp_0|^2}{\hbar}
\right),
\end{equation}
with
$\bx_0=(1.5,-1.0,0.5)^\top$
and
$\bp_0=(0.5,0.25,-0.75)^\top$.
The remaining parameters are
\begin{equation}
m=1,\qquad
\hbar=0.5,\qquad
\gamma=0.1,\qquad
D_{pp}=D_{xx}=0.08,\qquad
D_{xp}=0.
\end{equation}
These coefficients satisfy the Lindblad condition \cref{eq_lindblad_condition} so that the corresponding continuous density-operator evolution is completely positive. Thus, this example provides a genuinely three-dimensional open-quantum-system test of the proposed low-rank method.

The computational domain is truncated to
\[
[x_{\min},x_{\max}]^3\times[p_{\min},p_{\max}]^3
=
[-6,6]^3\times[-6,6]^3,
\]
with $N_x=N_p=64$ grid points in each direction.
We perform the simulations using two CFL coefficients, $\mathrm{CFL}=0.05$ and $\mathrm{CFL}=0.025$.
The prescribed TT rank is $(24,24,24,24,24)$.
Under this discretization, a direct full-tensor representation would require
$64^6\approx6.87\times10^{10}$ stored values.
In contrast, the TT representation requires only
\[
2(64\times24)+4(64\times24^2)=150528
\]
stored values.
Thus, the tensor-train format reduces the storage requirement by approximately five orders of magnitude.
To the best of our knowledge, numerical simulations of the time-dependent Wigner--Fokker--Planck equation in a full three-dimensional position--momentum setting remain very limited, and this example is intended to demonstrate the feasibility of the proposed method in such a high-dimensional regime.

\begin{figure}[t]
    \centering
    \includegraphics[width=0.8\linewidth]{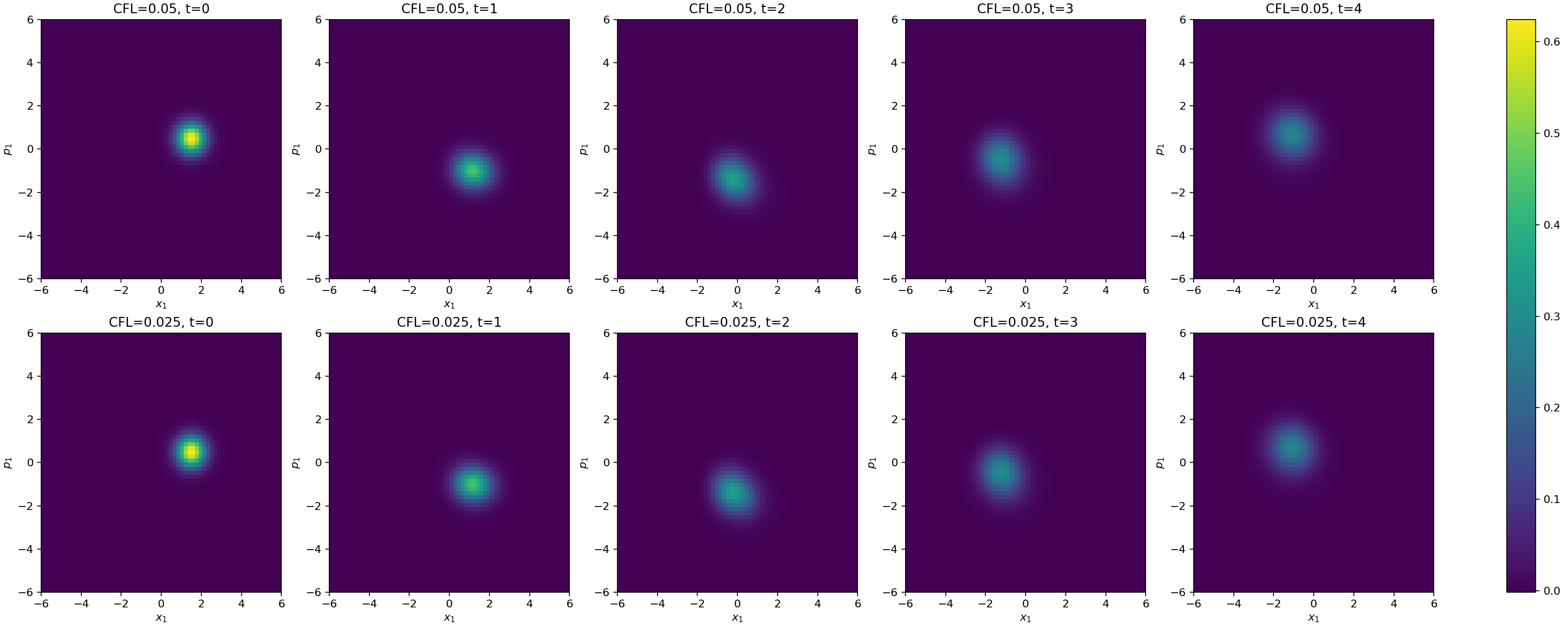}
    \caption{Evolution of $W_1(t,x_1,p_1)$ for different CFL coefficients.}
    \label{fig_example_3D3V_W1}
\end{figure}

\begin{figure}[t]
    \centering
    \includegraphics[width=0.8\linewidth]{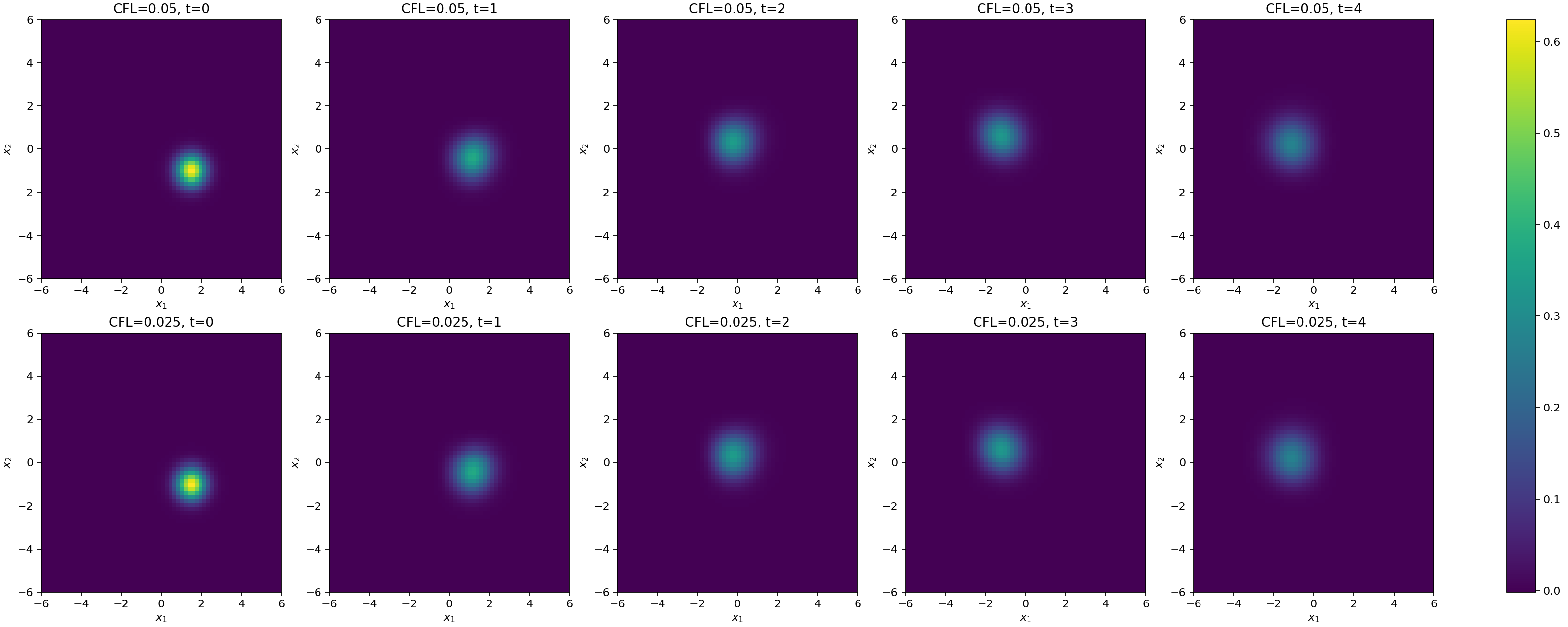}
    \caption{Evolution of $P_{12}(t,x_1,x_2)$ for different CFL coefficients.}
    \label{fig_example_3D3V_P12}
\end{figure}

Although the full state evolves in a six-dimensional phase space, physical information is often extracted through reduced phase-space distributions and spatial marginals.
Such quantities can be evaluated efficiently from the TT representation without reconstructing the full six-dimensional tensor.
We would like to examine two reduced quantities obtained by integrating out selected phase-space variables.
The reduced Wigner function associated with $(x_1,p_1)$ is defined by
\begin{equation}
    W_1(t,x_1,p_1)
    =
    \int_{\mathbb{R}^4}
    W(t,\bx,\bp)\,
    \dd x_2\,\dd x_3\,\dd p_2\,\dd p_3,
\end{equation}
while the reduced spatial density in $(x_1,x_2)$ is
\begin{equation}
    P_{12}(t,x_1,x_2)
    =
    \int_{\mathbb{R}^4}
    W(t,\bx,\bp)\,
    \dd x_3\,\dd\bp.
\end{equation}
Their evolutions are shown in
\cref{fig_example_3D3V_W1,fig_example_3D3V_P12}.
The numerical solutions obtained with the two CFL coefficients exhibit very similar evolution, indicating consistency with respect to the time-step refinement.

We next examine the evolution of the total mass,
\begin{equation}
M(t)
=
\int_{\mathbb{R}^3}
\int_{\mathbb{R}^3}
W(t,\bx,\bp)\,\dd\bx\,\dd\bp.
\end{equation}
For the Wigner--Fokker--Planck equation considered here, the total mass is conserved analytically.
The present numerical method, however, does not impose exact mass conservation.
Monitoring $M(t)$ therefore provides a useful diagnostic for assessing the accumulated numerical error.
Although the physical Wigner function and its total mass are real-valued, the Fourier representation of the Wigner potential and the corresponding low-rank contractions involve complex-valued intermediate quantities.
We therefore monitor both the real and imaginary parts of the numerically computed mass.
Since the initial Wigner function is normalized, $M(0)=1$.

\begin{figure}[t]
\centering
\includegraphics[width=0.8\linewidth]{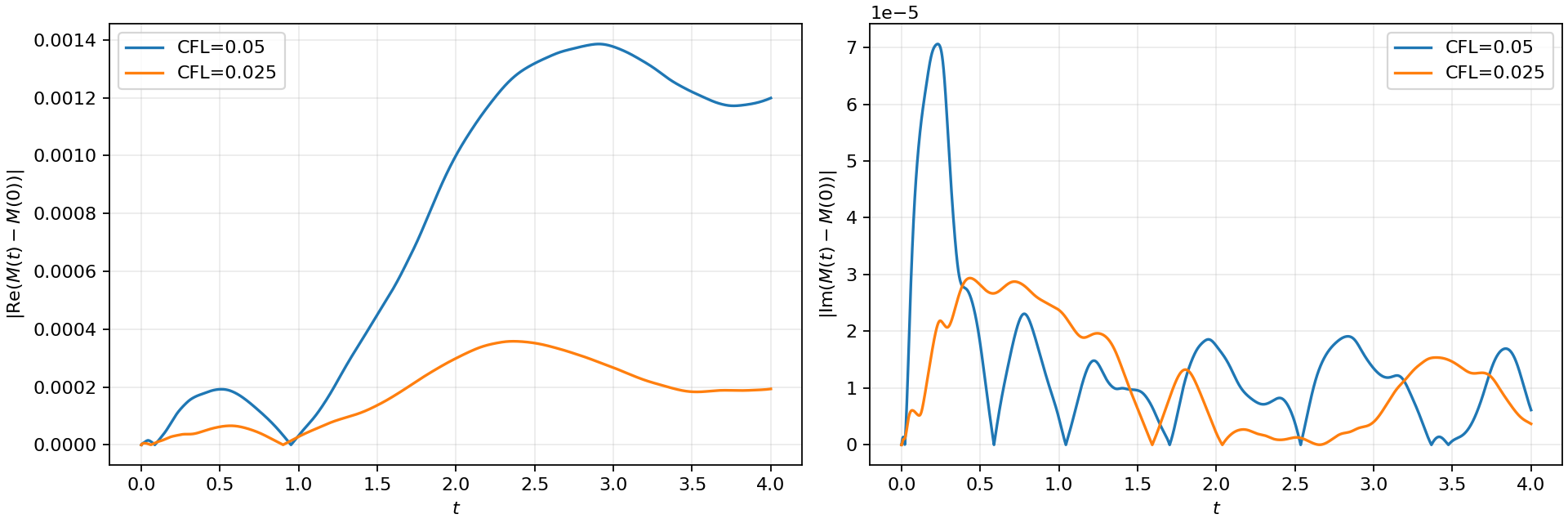}
\caption{Evolution of the absolute values of the real and imaginary parts of the mass error $M(t)-M(0)$ for different CFL coefficients.}
\label{fig_example_3D3V_mass}
\end{figure}

The evolution of the mass error is shown in \cref{fig_example_3D3V_mass}.
Although exact conservation is not enforced, the deviation of the real part of $M(t)$ from unity remains small throughout the simulation, while the imaginary part remains close to zero.
Moreover, the simulation with the smaller CFL coefficient exhibits a smaller mass error, consistent with the reduction of the time-discretization error as the time step decreases.

Unlike low-rank schemes formulated directly in the density-matrix setting that enforce complete positivity and trace preservation by construction \cite{appelo2025kraus,delmastro2026completely}, the present Wigner-space projector-splitting method does not impose these properties at the discrete level.
Nevertheless, the small mass errors observed here indicate that the trace normalization is accurately maintained in this experiment.
Developing a conservative tensor-train formulation that preserves the total mass, and possibly additional physical moments, remains an important direction for future work.

In addition, we vary the grid resolution and the prescribed TT rank to examine the empirical computational scaling of the proposed method.
For simplicity, we take $N_x=N_p=N$ in all tests.
\Cref{table_wall_clock_time} reports the wall-clock time for a single time-step propagation, obtained by averaging the computational time over 20 consecutive time steps.
To estimate the dependence on $N$ and $r$, we fit the measured timings to a power-law model of the form $CN^\alpha r^\beta$ using least squares.
The resulting fit gives an empirical scaling of approximately
$\mathcal{O}(N^{1.07}r^{2.12})$.
This behavior is consistent with the complexity estimate in \cref{eq_computational_cost}.
In particular, the nearly linear dependence on $N$ and approximately quadratic dependence on $r$ suggest that, over the range of parameters considered here, a significant portion of the computational cost is associated with the
$\mathcal{O}(d^2\varrho r^2N_p\log N_p)$
term arising from the FFT-related operations in the construction of the potential environments.
The weak logarithmic dependence on $N_p$ is difficult to distinguish from linear scaling over the present range of grid sizes, which is consistent with the observed exponent $1.07$.
The observed exponent $2.12$ in $r$ is also close to the quadratic rank dependence of this contribution.

\begin{table}[t]
\centering
\begin{tabular}{c|ccccccccc}
\hline
$N$ & $r=24$ & $r=26$ & $r=28$ & $r=30$ & $r=32$ & $r=34$ & $r=36$ & $r=38$ & $r=40$ \\
\hline
64
& 0.926 & 1.041 & 1.164 & 1.354 & 1.604 & 1.830 & 2.008 & 2.272 & 2.504 \\
128
& 1.761 & 2.103 & 2.387 & 2.821 & 3.369 & 3.873 & 4.292 & 4.855 & 5.383 \\
256
& 3.877 & 4.613 & 5.223 & 6.042 & 7.251 & 8.056 & 8.823 & 9.977 & 10.954 \\
512
& 7.940 & 9.265 & 10.814 & 12.728 & 15.471 & 16.986 & 19.011 & 21.219 & 23.722 \\
\hline
\end{tabular}
\caption{Average wall-clock time (in seconds) for one time-step propagation with different grid resolutions $N_x=N_p=N$ and prescribed TT ranks $r$.}
\label{table_wall_clock_time}
\end{table}

\section{Conclusion}
\label{sec_conclusion}
In this work, we develop a dynamical tensor train method for the numerical simulation of the Wigner--Fokker--Planck equation in high-dimensional phase space.
The Wigner function is represented and evolved directly in TT form using the projector-splitting integrator.
For the nonlocal Wigner pseudo-differential operator, we introduce low-rank CP and TT representations of the potential difference and incorporate the resulting projected operator via a matrix-free exponential action.
The other terms in the WFP equations are evaluated and projected directly in the phase-space variables.
The method avoids constructing the full phase-space tensor throughout the computation and substantially reduces memory and computational costs.
Numerical experiments in one-, two- and three-dimensional settings are included to demonstrate the accuracy and efficiency of the method. 

Several directions remain for future investigation.
The present method uses fixed TT ranks and a first-order projector-splitting time discretization.
Rank-adaptive frameworks and higher-order time integrators could further improve its flexibility and efficiency.
In addition, the current low-rank discretization does not enforce trace preservation or complete positivity at the discrete level, although these properties are satisfied by the underlying continuous dynamics under the corresponding Lindblad conditions.
Developing structure-preserving tensor-train methods that retain these fundamental properties of open quantum evolution is therefore an important direction for future work.
Another interesting direction is to investigate alternative low-rank tensor formats, such as Tucker and hierarchical Tucker representations, as well as other dynamical low-rank integrators, such as the basis-update-and-Galerkin (BUG) integrator, for the Wigner--Fokker--Planck equation.
A systematic comparison of these representations and integrators would help clarify their relative advantages in terms of accuracy, efficiency, rank growth, and structure preservation.

\section*{Acknowledgement}
This work was partially supported by DOE grant DE-SC0023164, NSF grants DMS-2409858 and IIS-2433957, and AFOSR grants FA9550-24-1-0254 and FA9550-26-1-B161.

\bibliographystyle{abbrv}
\bibliography{myBib}

\end{document}